\documentstyle[12pt,psfig]{article}

\def\hatM{{\widehat{M}}}
\def\hatv{{\widehat{v}}}
\def\hattt{{\widehat{\calt}}}
\def\caltbar{\overline{\calt}}
\def\hatpp{{\widehat{\pp}}}
\def\interior{{\rm Int}}
\def\index{{\rm ind}}
\def\vbar{{\overline v}}
\def\Wbar{{\overline W}}
\def\Kbar{{\overline K}}

\def\hbar{{\overline h}}
\def\GL{{\rm GL}}

\def\Ccell{{\rm C}^{\rm cell}}
\def\Cphi{{\rm C}^\varphi}

\def\basis{{\hbox{\Got b}}}
\def\hbasis{{\hbox{\Got h}}}
\def\gbasis{{\hbox{\Got g}}}

\def\ristr#1{\big|_{#1}}

\def\stwotriv{S^2_{\rm triv}}

\reversemarginpar

\makeatletter
\def\@begintheorem#1#2{\it \trivlist \item[\hskip \labelsep{\bf #1\ #2.}]}
\newtheorem{teo}{Theorem}[section]
\newtheorem{rem}[teo]{Remark}
\newtheorem{lem}[teo]{Lemma}

\newtheorem{cor}[teo]{Corollary}

\newtheorem{prop}[teo]{Proposition}

\def\finedim#1{{\hfill\hbox{\enspace\fbox{\ref{#1}}}}\vspace{5pt}}
\def\dim#1{\vspace{1pt}\noindent{\it Proof of} {\hspace{2pt}}\ref{#1}.}
\def\compo{\,{\scriptstyle\circ}\,}
\def\cont{{\rm C}}

\font\scpicc=cmcsc10
\font\sc=cmcsc10 scaled 1200

\newfont{\Bbb}{msbm10 scaled 1200}
\def\mr{{\hbox{\Bbb R}}}

\def\mz{{\hbox{\Bbb Z}}}

\newfont{\mycal}{eusm10 scaled 1200}

\newfont{\Got}{eufm10 scaled 1200}

\def\cc{{\cal C}}

\def\uu{{\cal U}}
\def\ss{{\cal S}}

\def\calt{{\cal T}}

\def\hh{{\cal H}}
\def\pp{{\cal P}}

\def\dd{{\cal D}}

\def\eul{{\rm Eul}}
\def\euls{{\rm Eul}^{\rm s}}
\def\eulc{{\rm Eul}^{\rm c}}
\def\Thetas{\Theta^{\rm s}}
\def\Thetac{\Theta^{\rm c}}
\def\alphas{\alpha^{\rm s}}
\def\alphac{\alpha^{\rm c}}

\title{Reidemeister Torsion of 3-Dimensional\\ Euler Structures
with Simple Boundary Tangency\\ and Pseudo-Legendrian Knots}

\author{Riccardo Benedetti\qquad Carlo Petronio\thanks{The second named
author gratefully acknowledges financial support by GNSAGA-CNR}}

\begin{document}

\maketitle

\noindent{\small{\scpicc Abstract}. We generalize Turaev's definition of torsion invariants of pairs $(M,\xi)$, where $M$ is a 3-dimensional manifold and $\xi$ is an Euler structure on $M$ 
(a non-singular
vector field up to homotopy relative to $\partial M$  and local modifications in $\interior(M)$). Namely, we
allow $M$ to have arbitrary boundary and $\xi$ to have simple (convex and/or concave) tangency circles to
the boundary.  We prove that Turaev's $H_1(M)$-equivariance formula holds also in  our generalized context.
Our torsions apply in particular to (the exterior of) pseudo-Legendrian knots 
({\em i.e.}~knots transversal to a given vector field),
and hence to Legendrian knots in contact
3-manifolds. We show that torsion, as an
absolute invariant, contains a lifting to pseudo-Legendrian knots of the classical
Alexander invariant. We also precisely analyze the information carried by torsion
as a relative invariant of pseudo-Legendrian knots which are framed-isotopic.
Using branched standard spines to describe
vector fields we
show how to explicitly invert Turaev's reconstruction map  from combinatorial to smooth Euler structures,
thus making the computation of torsions a more effective one. As an example we work out a specific
calculation.}

\vspace{.5cm}


\vspace{.5cm}

\section*{Introduction}
Reidemeister torsion is a classical yet very vital topic in
3-dimensional topology, and it was recently used in a variety of important
developments. To mention a few, torsion is a fundamental ingredient of the
Casson-Walker-Lescop invariants (see {\em e.g.}~\cite{lescop}), and more
generally of the perturbative approach to quantum invariants (see {\em
e.g.}~\cite{lemuoh}). Relations have been pointed out between torsion and
hyperbolic geometry~\cite{porti}. Turaev's torsion of Euler 
structures~\cite{turaev:Euler} has recently been 
recognized by Turaev himself (\cite{turaev:spinc},~\cite{turaev:nuovo}) to have deep
connections with the Seiberg-Witten
invariants of ${\rm Spin}^{\rm c}$-structures on 3-manifolds, after 
the proof of Meng and Taubes~\cite{meng:taub} that a suitable combination of these
invariants can be identified with the classical Milnor torsion.

Turaev's theory~\cite{turaev:Euler} actually exists in all dimensions. We
quickly review it before proceeding. A {\em smooth Euler structure} $\xi$ on a
compact oriented manifold $M$, possibly with $\partial
M=\emptyset$, is a non-singular vector field on $M$ viewed up to local
modifications  in $\interior(M)$ and  homotopy relative to $\partial M$. 
Orientability of $M$ is not strictly necessary, but we find
it convenient to assume it. Turaev
allows only ``monochromatic'' boundary components, {\it i.e.}~black ones (on
which the field points outwards) and white ones (on which it points inwards).
This implies the constraint that $\chi(M,W)=0$, where $W$ is the white portion
of  $\partial M$, but in~\cite{turaev:spinc} and~\cite{turaev:nuovo} Turaev
only focuses on the more specialized case where $M$ is 3-dimensional and closed
or bounded by tori. In all dimensions, the set $\euls(M,W)$ of smooth Euler
structures compatible with $(M,W)$ is an affine space over  $H_1(M;\mz)$. The
two main ingredients of Turaev's theory are as follows. First, he defines a
certain set of $1$-chains, called the space $\eulc(M,W)$ of {\em combinatorial}
Euler structures compatible with $(M,W)$, he shows that this is again affine
over $H_1(M;\mz)$,  and he describes an $H_1(M;\mz)$-equivariant bijection
$\Psi:\eulc(M,W)\to\euls(M,W)$ called the {\em reconstruction map}. Second, for
$\xi\in\eulc(M,W)$ and for  any representation $\varphi$ of $\pi_1(M)$ into the
units of a suitable ring $\Lambda$ he defines a torsion  invariant
$\tau^\varphi(M,\xi)$, or more generally $\tau^\varphi(M,\xi,\hbasis)$, with
values in $K_1(\Lambda)/(\pm1)$. This invariant is by definition a lifting of
the classical Reidemeister torsion (see~\cite{milnor}) $\tau^\varphi(M)\in
K_1(\Lambda)/(\pm\varphi(\pi_1(M)))$, and it satisfies the
$H_1(M;\mz)$-equivariance formula \begin{equation} \label{tau:equivariance}
\tau^\varphi(M,\xi',\hbasis)=\tau^\varphi(M,\xi,\hbasis)\cdot \varphi(\xi'-\xi)
\end{equation} where $\xi'-\xi\in H_1(M;\mz)$. For $\xi\in\euls(M,W)$ one
defines $\tau^\varphi(M,\xi)$ as $\tau^\varphi(M,\Psi^{-1}(\xi))$, and the 
$H_1(M;\mz)$-equivariance of the reconstruction map $\Psi$ implies that 
formula~(\ref{tau:equivariance}) holds also for smooth structures. We emphasize
that the definition of $\Psi$ is based on an explicit geometric construction,
but its bijectivity is only established through $H_1(M;\mz)$-equivariance. This
makes the definition of torsion for smooth structures somewhat implicit.

In the present paper, and in other papers in preparation, we are concerned with
generalizations and improvements of Turaev's theory. Here we consider
3-manifolds. This work had two main initial aims. Our first aim was to find a
geometric description of the map $\Psi^{-1}$, and hence to turn the computation
of Turaev's torsion into a more effective procedure, using our 
encoding~\cite{lnm} of non-singular vector fields up to homotopy (also
called ``combings'') in terms of branched standard spines. Our second aim was
to define torsion invariants of pseudo-Legendrian pairs $(v,L)$, 
consisting of a link $L$ transversal to a non-singular vector field $v$, viewed up to 
{\em pseudo-Legendrian isotopy}, namely transversality-preserving
simultaneous isotopy of $L$ and homotopy of $v$. For a given $v$ we will just say
that $L$ is pseudo-Legendrian in $(M,v)$, or just in $v$. Note that
an ordinary Legendrian link in a given oriented contact structure $\xi$
is pseudo-Legendrian in $\xi^\perp$, and Legendrian isotopy implies pseudo-Legendrian
isotopy.
A specific motivation to look for invariants of pseudo-Legendrian
links comes  from the remarkable relation recently discovered by Fintushel and
Stern~\cite{fist} between  the Alexander polynomial ({\em i.e.}~Milnor torsion)
of a knot $K\subset S^3$ and  (a suitable combination of) the Seiberg-Witten
invariants of the  ``surgered'' $4$-manifold $X_K$ obtained using $K$ (and a
suitable base $4$-manifold $X$). Both our initial aims lead us to consider
Euler structures on $3$-manifolds $M$ (without restrictions on $\partial M$)
allowing simple tangency circles to $\partial M$ of  {\em concave} type (see
Fig.~\ref{conc:conv:tang} below). On the other  hand it turns out that, to
define torsion, the natural objects to deal with  are Euler structures with
{\em convex} tangency  circles. It is a fortunate fact, peculiar of dimension
3, that there is a canonical way to associate  a convex field to any {\em
simple} ({\em i.e.}~mixed concave and convex) one. This allows to define
torsion for all smooth simple Euler structures, and eventually to achieve both
the objectives we had in mind.

Let us now summarize the contents of this paper. The foundational part of
our work consists in extending to the
context of Euler structures with simple tangency the notions of combinatorial
structure $\eulc$ and reconstruction map $\Psi$. 
This part follows the same scheme as~\cite{turaev:Euler}
and relies on technical results of Turaev. Our main contribution here
is the proof that the natural transformations of a concave structure
into a convex one, viewed at the smooth level and at the combinatorial level,
actually correspond to each other under the reconstruction map
(Theorem~\ref{diagram:commutes}).
After setting the foundations, we prove the following main result
(stated informally here:
see Sections~\ref{Eul:def:section} 
and~\ref{spines:section} for precise definitions and statements).

\begin{teo}\label{informal:algorithmic}
Let $\xi$ be an Euler structure with concave
tangency circles. If $P$ is a branched standard spine 
which represents $\xi$, then $P$ allows to explicitly 
find a representative of $\Psi^{-1}(\xi)\in\eulc$, and hence to compute
the torsion of $\xi$ in terms of the finite combinatorial data which encode $P$.
\end{teo}

After the foundations, we concentrate on pseudo-Legendrian knots $(v,K)$, assuming
for simplicity the ambient manifold to be closed. The connection
comes from the fact that the restriction of $v$ defines
a concave Euler structure on the exterior $E(K)$ of $K$, with two parallel tangency lines
on $\partial E(K)$ determined by the framing defined by $v$ on $K$. We 
show that torsion, as an {\em absolute} invariant, contains (in a suitable sense)
a lifting to pseudo-Legendrian knots of the classical Alexander invariant.
Then we carefully analyze the {\em relative} information carried by torsion
for pseudo-Legendrian pairs $(v_0,K_0)$ and $(v_1,K_1)$ such that $(v_0,v_1)$ 
are homotopic to each other and $(K_0,K_1)$ are framed-isotopic to each other.
A relevant point which emerges from this analysis is that in general torsion does not
provide a single-valued relative invariant, because the
action of a certain mapping class group (which depends on the framed isotopy class only)
must be taken into account. This leads us to the notion of `good' framed knots,
for which the action is trivial, and the study of torsion is simpler. 
We show that many knots are good (for instance, all knots in a homology sphere are good,
and most knots with hyperbolic complement are good).
Concentrating on good knots we then prove that in a homology sphere the relative torsion
of two knots essentially coincides with the difference of their rotation numbers
(Maslov indices), so torsion basically detects whether the knots are isotopic through
pseudo-Legendrian immersions.

Moreover we analyze the effect on torsion
of the framed first Reidemeister move (which does not change the framing
but locally changes the winding number by $\pm2$), and we show that for
homology spheres the winding number is just the difference of Maslov indices,
thus getting an alternative proof of the relation between torsion and rotation number.
Using the fact (proved in~\cite{second:paper}) that framed isotopy is
generated by pseudo-Legendrian isotopy and the framed first Reidemeister move,
we then obtain several interesting consequences, among which we state the following:

\begin{teo}\label{informal:good:for:knots}
Consider pseudo-Legendrian knots
$(v_0,K_0)$ and $(v_1,K_1)$. Assume that $K_0$ is good and that the meridian of $K_0$ has
infinite order in $H_1(E(K_0);\mz)$. Then the knots are pseudo-Legendrian
isotopic if and only if they have trivial relative torsion invariants.
\end{teo}

This paper is organized
as follows. In Section~\ref{Eul:def:section} we provide the
formal definitions of smooth and combinatorial Euler structure.
In Section~\ref{torsion:def:section} we introduce torsion and
state the equivariance property. In Section~\ref{spines:section} we
show how branched standard spines can be used for computing torsion.
In Section~\ref{knots:section} we specialize to pseudo-Legendrian 
knot exteriors and analyze torsion both as an absolute and as a relative invariant.
In Section~\ref{exa:section} we carry out a specific computation using
the technology of Section~\ref{spines:section}. 
In Sections~\ref{Eul:def:section}
to~\ref{knots:section} proofs which are long and require the introduction of 
ideas and techniques not used elsewhere are omitted. Section~\ref{proofs} contains all
these proofs.

We conclude this introduction by announcing related results which we have
recently obtained and partially written down. In~\cite{second:paper} we extend
to the case with boundary our combinatorial presentation~\cite{lnm} of combed
manifolds in terms of branched spines, and we provide similar presentations
of framed and pseudo-Legendrian links, using $\cont^1$ diagrams on branched 
spines. Some results from~\cite{second:paper} are actually used also in the
present paper (see Sections~\ref{spines:section} and~\ref{knots:section}). 
In~\cite{third:paper} we use the results of~\cite{second:paper} to develop
an approach to torsion entirely based on combinatorial techniques, getting
slightly different generalizations of Turaev's theory. In~\cite{fourth:paper} we
generalize the theory of Euler structures and (with some restrictions) of
torsion to all dimensions and allowing any generic (Whitney-Morin-type)
tangency to the boundary.

\section{Euler structures}\label{Eul:def:section}
In this section we define smooth and combinatorial
Euler structures and explain their correspondence.
Fix once and for ever a compact oriented 3-manifold $M$, possibly with 
$\partial M=\emptyset$.
Using the {\it Hauptvermutung}, we will always freely intermingle the
differentiable, piecewise linear and topological viewpoints. Homeomorphisms
will always respect orientations. All vector fields mentioned in this paper will be 
non-singular, and they will be termed just {\em fields} for the sake of brevity.

\begin{figure}
\centerline{\psfig{file=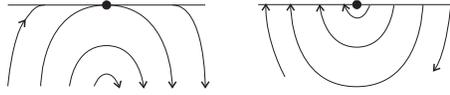,width=6cm}}
\caption{\label{conc:conv:tang} Convex (left) and concave (right) tangency to the boundary.}
\end{figure}

\paragraph{Smooth and combinatorial Euler structures}
We will call {\em boundary pattern} on $M$ a partition $\pp=(W,B,V,C)$ of
$\partial M$ where $V$ and $C$ are finite unions of 
disjoint circles, and $\partial W=\partial B=V\cup C$. In particular, $W$ and $B$ are
interiors of compact surfaces embedded in $\partial M$. Even if $\pp$ can actually 
be determined by less data, {e.g.} the pair $(W,V)$,
we will find it convenient to refer to $\pp$ as a quadruple.
Points of $W$, $B$, $V$ and $C$ will be called {\em white}, {\em black}, 
{\em convex} and {\em concave} respectively.
We define the set
of {\em smooth Euler structures} on $M$ compatible with $\pp$, denoted by
$\euls(M,\pp)$, as the set of equivalence classes of
fields on $M$ which point inside on $W$, point
outside on $B$ and have simple tangency to $\partial M$ 
of {\em convex} type along $V$ and {\em concave} type 
along $C$,
as shown in a cross-section in Fig.~\ref{conc:conv:tang}.
Two such fields are equivalent if they are obtained from each other by 
homotopy through fields of the same type and modifications supported into
interior balls. The following
variation on the Poincar\'e-Hopf formula is established in Section~\ref{proofs}:

\begin{prop}\label{p:h:formula}
$\euls(M,\pp)$ is non-empty if and only if $\chi(\Wbar)=\chi(M)$. 
\end{prop}

\noindent We remark here that $\chi(\Wbar)=\chi(W)$, $\chi(\overline{B})=\chi(B)$, $\chi(V)=\chi(C)=0$ and
$\chi(W)+\chi(B)=\chi(\partial M)=2\chi(M)$, so there are various ways to
rewrite the relation $\chi(\Wbar)=\chi(M)$, the most intrinsic of which is actually
$\chi(M)-(\chi(\Wbar)-\chi(C))=0$ (see below for the reason).

Now, given $\xi,\xi'\in\euls(M,\pp)$ we can choose generic representatives $v,v'$, so
that the set of points of $M$ where $v'=-v$ is a union of loops contained in the
interior of $M$. A standard procedure allows to give these loops a canonical orientation,
thus getting an element $\alphas(\xi,\xi')\in H_1(M;\mz)$. The following result is easily
obtained along the lines of the well-known analogue for closed manifolds.

\begin{lem}
$\alphas$ is well-defined and turns $\euls(M,\pp)$ into an affine space over $H_1(M;\mz)$.
\end{lem}

A (finite) cellularization $\cc$ of $M$ 
is called {\em suited} to $\pp$ if $V\cup C$ is
a subcomplex, so $W$ and $B$ are unions of cells. 
Here and in the sequel by ``cell'' we will always mean an {\em open} one.
Let such a $\cc$ be given.
For $\sigma\in\cc$ define $\index(\sigma)=(-1)^{{\rm dim}(\sigma)}$. We define
$\eulc(M,\pp)_\cc$ as the set of equivalence classes of integer singular 1-chains $z$ in $M$
such that 
$$\partial z=\sum_{\sigma\subset M\setminus (W\cup V)}\index(\sigma)\cdot p_\sigma$$
where $p_\sigma\in\sigma$ for all $\sigma$. Two chains $z$ and $z'$ 
with $\partial z=\sum\index(\sigma)\cdot p_\sigma$ and 
$\partial z'=\sum\index(\sigma)\cdot p'_\sigma$
are defined to be equivalent if there exist $\delta_\sigma:([0,1],0,1)\to
(\sigma,p_\sigma,p'_\sigma)$ such that
$$z-z'+\sum_{\sigma\subset M\setminus (W\cup V)}\index(\sigma)\cdot \delta_\sigma$$
represents $0$ in $H_1(M;\mz)$. Elements of
$\eulc(M,\pp)_\cc$ are called {\em combinatorial Euler structures} relative to $\pp$ and $\cc$,
and their representatives are called {\em Euler chains}.
The definition implies that, for $\xi,\xi'\in\eulc(M,\pp)_\cc$, their difference $\xi-\xi'$
can be defined as an element $\alphac(\xi,\xi')$ of $H_1(M;\mz)$. The following is easy:

\begin{lem}
$\eulc(M,\pp)_\cc$ is non-empty if and only if $\chi(\Wbar)=\chi(M)$, and in this case
$\alphac$ turns it into an affine space over $H_1(M;\mz)$.
\end{lem}

\noindent Since $\Wbar=W\cup V\cup C$, the alternating sum of dimensions of cells in $W\cup V$ is 
intrinsically interpreted as $\chi(\Wbar)-\chi(C)$, which explains why the most meaningful
way to write the relation $\chi(\Wbar)=\chi(M)$ is $\chi(M)-(\chi(\Wbar)-\chi(C))=0$.
From now on we will always assume that this relation holds.
Turaev~\cite{turaev:Euler} only considers the case where $V=C=\emptyset$, so 
$W=\overline{W}$ and $B=\overline{B}$, and our relation takes the usual form
$\chi(M,W)=0$.
The following result was established by Turaev in~\cite{turaev:Euler}
in his setting, but the proof
extends {\em verbatim} to our context, so we omit it.
Only the first assertion is hard.
We state the other two because we will use them.

\begin{prop}\label{combin:prop:statement}
\begin{enumerate}
\item If $\cc'$ is a subdivision of $\cc$ then there exists a canonical
$H_1(M;\mz)$-isomorphism $\eulc(M,\pp)_\cc\to\eulc(M,\pp)_{\cc'}$. In particular
$\eulc(M;\mz)$ is canonically defined up to $H_1(M;\mz)$-isomorphism independently 
of the cellularization.
\item\label{connected:spider:point} If 
$\cc$ is a cellularization of $M$ suited to $\pp$
and $x_0\in M$ is an assigned point, any element of $\eulc(M,\pp)$ can be represented,
with respect to $\cc$, as a sum $\sum_{\sigma\subset M\setminus(W\cup V)}
\index(\sigma)\cdot \beta_\sigma$
with $\beta_\sigma:([0,1],0,1)\to(M,x_0,\sigma)$.
\item\label{bary:sub:point} If $\calt$ is a triangulation of $M$ suited to $\pp$,
any element of $\eulc(M,\pp)$ can be represented, with respect to $\calt$, as 
a simplicial $1$-chain in the first barycentric subdivision of $\calt$.
\end{enumerate}
\end{prop}

\noindent Our first main result, proved in Section~\ref{proofs},
is the extension to the case under consideration of
Turaev's correspondence between $\eulc$ and $\euls$.

\begin{teo}\label{reconstruction:statement}
There exists a canonical $H_1(M;\mz)$-equivariant isomorphism 
$$\Psi:\eulc(M,\pp)\to\euls(M,\pp).$$
\end{teo}

\noindent The definition of $\Psi$ is based on an explicit geometric construction,
but its bijectivity is only established through $H_1(M;\mz)$-equivariance.
As already mentioned in the introduction, this makes in general a very difficult task to
determine the inverse of $\Psi$. One of the features of this paper is the description
of $\Psi^{-1}$ in terms of the combinatorial encoding of fields by 
means of branched spines: Theorem~\ref{spider:structure:teo} 
describes $\Psi^{-1}$ when $\pp$ is concave,
and Theorem~\ref{diagram:commutes} shows that from a general $\pp$ we can
effectively pass to a unique convex $\pp$, and hence to a unique concave $\pp$,
and conversely.

In view of Theorem~\ref{reconstruction:statement}, when no confusion risks to arise,
we shortly write $\eul(M,\pp)$ for either $\euls(M,\pp)$ or $\eulc(M,\pp)$,
and $\alpha$ for the map giving the affine $H_1(M;\mz)$-structure on this space.

\begin{figure}
\centerline{\psfig{file=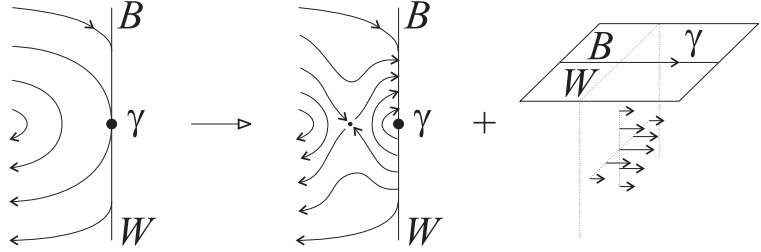,width=10cm}}
\caption{\label{conc:to:conv} Turning a concave tangency circle $\gamma$ into a convex one: the 
apparent singularity in the cross-section is removed by adding a small bell-shaped field 
directed parallel to $\gamma$, {\em i.e.}~orthogonal to the cross-section.}
\end{figure}

\paragraph{Convex Euler structure associated to an arbitrary one}
Let $M$ and $\pp=(W,B,V,C)$ be as in the definition of $\eul(M,\pp)$.
The pattern $\theta(\pp)=(W,B,V\cup C,\emptyset)$ is a convex one canonically
associated to $\pp$. We define a map 
$$\Thetas:\euls(M,\pp)\to\euls(M,\theta(\pp))$$
as geometrically described in Fig.~\ref{conc:to:conv}. Concerning
this figure, note that the loops in $C$ can be oriented as 
components of the boundary of $B$, which is oriented as a subset of the boundary of $M$.

\begin{lem}\label{concave:convex:smooth}
$\Thetas$ is a well-defined $H_1(M;\mz)$-equivariant bijection.
\end{lem}

\dim{concave:convex:smooth}
The first two properties are easy and imply the third property. The inverse of $\Thetas$
may actually be described geometrically by a figure similar to Fig.~\ref{conc:to:conv}, but we leave this to the reader.
\finedim{concave:convex:smooth}
 
We define now a combinatorial version of $\Thetas$. Consider a cellularization $\cc$ suited to 
$\pp$, and denote by $\gamma_1,\dots,\gamma_n$ the 1-cells contained in $C$. We choose the 
parameterizations $\gamma_j:(0,1)\to C$ so that they respect the natural orientation of $C$ 
already discussed above, and we extend the $\gamma_j$ to $[0,1]$, without changing notation.
Now let $z$ be an Euler chain relative to $\pp$. It 
easily seen that $z-\sum_{j=1}^n\gamma_j\ristr{[1/2,1]}$ is an Euler chain
relative to $\theta(\pp)$. Setting 
$$\Thetac([z])=\left[z-\sum_{j=1}^n\gamma_j\ristr{[1/2,1]}\right]$$ 
we get a map $\Thetac:\eulc(M,\pp)\to\eulc(M,\theta(\pp))$.

\begin{lem}\label{concave:convex:comb}
$\Thetac$ is a well-defined $H_1(M;\mz)$-equivariant bijection.
\end{lem}

\dim{concave:convex:comb}
Again, the first two properties are easy and imply the third one. 
\finedim{concave:convex:comb}

In Section~\ref{proofs} we will see the following:

\begin{teo}\label{diagram:commutes}
If $\Psi$ is the reconstruction map of Theorem~\ref{reconstruction:statement} then the
following diagram is commutative:
$$\matrix{
\eulc(M,\pp) & \stackrel{\Thetac}{\longrightarrow} & \eulc(M,\theta(\pp))\phantom{.} \cr
\Psi\downarrow\phantom{\Psi} & & \phantom{\Psi}\downarrow\Psi \cr
\euls(M,\pp) & \stackrel{\Thetas}{\longrightarrow} & \euls(M,\theta(\pp)). \cr}$$
\end{teo}

Using this result we will sometimes just write 
$\Theta:\eul(M,\pp)\to\eul(M,\theta(\pp))$.

\section{Torsion of an Euler structure}\label{torsion:def:section}
In this section we define torsion. We set up the usual
algebraic environment~\cite{milnor} in which 
torsion can be defined, fixing a ring $\Lambda$ with unit, with the property
that if $n$ and $m$ are distinct positive integers then $\Lambda^n$ and 
$\Lambda^m$ are not isomorphic as $\Lambda$-modules. The 
Whitehead group $K_1(\Lambda)$ is defined as the
Abelianization of $\GL_\infty(\Lambda)$, and $\Kbar_1(\Lambda)$ is the
quotient of $K_1(\Lambda)$ under the action of $-1\in\GL_1(\Lambda)=
\Lambda_*$. (Later in this paper the symbol $K_1$ will also be used for a knot,
but the meaning will always be clear from the context.)

We will directly define torsion only for a {\em convex} Euler structure,
but the definition easily extends to any Euler structure $\xi$ with simple
boundary tangency, taking the torsion of the convexified structure $\Theta(\xi)$.
So, we fix a manifold $M$, a {\em convex} boundary pattern $\pp=(W,B,V,\emptyset)$ on $M$,
a cellularization $\cc$ suited to $\pp$ and
a representation $\varphi:\pi_1(M)\to\Lambda_*$. 
We will denote by $\varphi$ again the extension
$\mz[\pi_1(M)]\to\Lambda$ (a ring homomorphism). 

We consider now the universal cover $q:\tilde{M}\to M$ and
the twisted chain complex $\Cphi_*(M,W\cup V)$, where
$\Cphi_i(M,W\cup V)$ is defined as 
$\Lambda\otimes_{\varphi}\Ccell_i(\tilde{M},q^{-1}(W\cup V);\mz)$,
and the boundary operator is induced from the ordinary boundary.
The homology of this complex is denoted by $H^\varphi_*(M,W\cup V)$ and
called the $\varphi$-twisted homology. We assume that each $H^\varphi_i(M,W\cup V)$
is a free $\Lambda$-module and fix a basis $\hbasis_i$.

\begin{rem}\label{twisted:remarks}{\em
\begin{enumerate}
\item\label{basepoint:remark} To have a formal completely intrinsic definition of  
$H^\varphi_*(M,W\cup V)$, one should fix from the
beginning a basepoint $x_0\in M$ for $\pi_1(M)$, and consider pointed
universal covers $q:(\tilde M,\tilde x_0)\to(M,x_0)$, because any two such covers
are {\em canonically} isomorphic, and the action of $\pi_1(M)$ on $\tilde{M}$ is
{\em canonically} defined on them. 
\item To define $H^\varphi_*(M,W\cup V)$
we have used in an essential way the fact that $W\cup V=\overline{W}$ is 
closed, because otherwise $\Cphi_*(M,W\cup V)$ cannot be defined.
\item\label{free:basis:remark} $\Cphi_i(M,W\cup V)$ is a free $\Lambda$-module, and
each $\mz[\pi_1(M)]$-basis of $\Ccell_i(\tilde{M},q^{-1}(W\cup V);\mz)$ determines a
$\Lambda$-basis of $\Cphi_i(M,W\cup V)$.
\item\label{varphibar:defin} If we compose $\varphi$ with the projection 
$\Lambda_*\to\Kbar_1(\Lambda)$ we get a homomorphism of $\pi_1(M)$ into an
{\em Abelian} group, so we get a homomorphism
$\overline{\varphi}:H_1(M;\mz)\to\Kbar_1(\Lambda).$
\end{enumerate}}
\end{rem}

\noindent Now let $\xi\in\eulc(M,\pp)$ and choose a representative of $\xi$ 
as in point~\ref{connected:spider:point} of 
Proposition~\ref{combin:prop:statement}, namely
$$\sum_{\sigma\in\cc,\ \sigma\subset M\setminus(W\cup V)}
\index(\sigma)\cdot \beta_\sigma$$
with $\beta_\sigma(0)=x_0$ for all $\sigma$, $x_0$ being a fixed point of $M$.
We choose $\tilde x_0\in q^{-1}(x_0)$ and consider the liftings
$\tilde\beta_\sigma$ which start at $\tilde x_0$. For $\sigma\subset M\setminus(W\cup V)$
we select its preimage $\tilde\sigma$ which contains $\tilde\beta_\sigma(1)$, and
define $\gbasis(\xi)$ as the collection of all these $\tilde\sigma$.
Arranging the $i$-dimensional elements of $\gbasis(\xi)$ in any order, by 
Remark~\ref{twisted:remarks}(\ref{free:basis:remark}) we
get a $\Lambda$-basis $\gbasis_i(\xi)$ of $\Cphi_i(M,W\cup V)$. 
We consider a set ${\tilde\hbasis}_i$ of elements of $\Cphi_i(M,W\cup V)$
which project to the fixed basis $\hbasis_i$ of $H^\varphi_i(M,W\cup V)$.

Now note that, given a free $\Lambda$-module $L$ 
and two finite bases $\basis=(b_k)$,
$\basis'=(b'_k)$ of $M$, the assumption made on $\Lambda$ guarantees that
$\basis $ and $\basis'$ have the same number of elements, so there
exists an invertible square matrix $(\lambda^h_k)$ such that
$b'_k=\sum_h\lambda^h_k b_h$. We will denote by
$[\basis'/\basis ]$ the image of $(\lambda^h_k)$ in $K_1(\Lambda)$.

\begin{prop}\label{defining:proposition}
If $\basis_i\subset\Cphi_i(M,W\cup V)$ is such that $\partial\basis_i$ is a
$\Lambda$-basis of $\partial(\Cphi_i(M,W\cup V))$, then
$(\partial\basis_{i+1})\cdot\tilde\hbasis_i\cdot\basis_i$ is a $\Lambda$-basis of
$\Cphi_i(M,W\cup V)$, and 
$$\tau^\varphi(M,\pp,\xi,\hbasis)=\pm
\prod_{i=0}^3\Big[ \Big( (\partial\basis_{i+1})\cdot\tilde\hbasis_i\cdot
\basis_i\Big) \;\Big/\;\gbasis_i(\xi)\Big]^{(-1)^i}\in
\Kbar_1(\Lambda)$$
is independent of all choices made. Moreover 
\begin{equation}
\tau^\varphi(M,\pp,\xi',\hbasis)=\tau^\varphi(M,\pp,\xi,\hbasis)\cdot 
\overline{\varphi}(\alphac(\xi',\xi)).\label{combinatorial:equivariance:formula}
\end{equation}
\end{prop}

\dim{defining:proposition}
The first assertion and independence of the $\basis_i$'s is purely algebraic and classical,
see~\cite{milnor}. Now note that $\xi\in\eulc(M,\pp)$ was used to select the bases 
$\gbasis_i(\xi)$. The $\gbasis_i(\xi)$ are of course not uniquely determined themselves,
but we can show that different choices lead to the same value of $\tau^\varphi$.

First of all, the arbitrary ordering in the $\gbasis_i(\xi)$ is inessential because
torsion is only regarded up to sign. Second, consider the effect of choosing a
different representative of $\xi$. This leads to a new family $\tilde\sigma'$ of
cells. If $\tilde\sigma'=a(\sigma)\cdot \tilde\sigma$, with $a(\sigma)\in\pi_1(M)$,
and $\overline{a}(\sigma)$ is the image in $H_1(M;\mz)$, we automatically have
$$\sum_{\sigma\subset M\setminus{W\cup V}}\index(\sigma)\cdot \overline{a}(\sigma)=0\in
H_1(M;\mz),$$
which allows to conclude that also the representative chosen is inessential.
The choice of the lifting $\tilde x_0$ can be shown to be inessential either
in the spirit of Remark~\ref{twisted:remarks}(\ref{basepoint:remark}),
or by showing that a simultaneous $a$-translation of all $\tilde\sigma$, for
$a\in\pi_1(M)$, multiplies the torsion by $\overline{\varphi}(a)^{\chi(M)-\chi(W\cup V)}=1$.

Formula~(\ref{combinatorial:equivariance:formula}) is readily established by choosing
representatives $\sum\index(\sigma)\cdot \beta_\sigma$ and 
$\sum\index(\sigma)\cdot \beta'_\sigma$ of $\xi$ and $\xi'$ 
such that $\beta'_\sigma=\beta_\sigma$ for all 
$\sigma$ but one.\finedim{defining:proposition}

Since the above construction uses the cellularization $\cc$ in a
way which may appear to be essential,
we add a subscript $\cc$ to the torsion we have defined.
The next result, which can be established following Turaev~\cite{turaev:Euler},
shows that dependence on $\cc$ is actually inessential.

\begin{prop}
Let $\cc$ and $\cc'$ be cellularizations suited to $\pp$. Assume that
$\cc'$ subdivides $\cc$, and consider the bijection
$\ss_{(\cc',\cc)}:\eulc(M,\pp)_\cc\to\eulc(M,\pp)_{\cc'}$ of
Proposition~\ref{combin:prop:statement}, and the canonical isomorphism
$j_{(\cc',\cc)}:H^\varphi_*(M,W\cup V)_\cc\to H^\varphi_*(M,W\cup V)_{\cc'}$. Then,
with obvious meaning of symbols we have:
$$\tau^\varphi_\cc(M,\pp,\xi,\hbasis)=
\tau^\varphi_{\cc'}(M,\pp,\ss_{(\cc',\cc)}(\xi),j_{(\cc',\cc)}(\hbasis)).$$
\end{prop}

It is maybe appropriate here to remark that the choice of a basis $\hbasis$  of
$H^\varphi_*(M,W\cup V)$ and the definition of $\tau^\varphi(M,\pp,\xi,\hbasis)$ implicitly 
assume a description of the universal cover of $M$, which is typically undoable in
practical cases. However,
if one starts from a representation of $\pi_1(M)$ into the units of a {\em commutative}
ring  $\Lambda$, {\em i.e.} a representation which factors through
 one of $H_1(M;\mz)$,
one can use from the very beginning the maximal Abelian rather
than the universal cover, which makes computations more feasible.

\begin{rem}{\em Turaev~\cite{turaev:Reidemeister} has shown that a homological
orientation yields a sign-refinement of torsion, {\em i.e.}~a lifting from 
$\Kbar_1(\Lambda)$ to $K_1(\Lambda)$. This refinement extends with minor
modifications to our setting of boundary tangency. This sign-refinement,
in the closed and monochromatic case, is often an essential component of the 
theory (for instance, it is crucial for the relation with
the 3-dimensional Seiberg-Witten invariants~\cite{turaev:spinc},~\cite{turaev:nuovo}
and for the definition of the Casson invariant~\cite{lescop}), but we will
not address it in the present paper.}
\end{rem}

\paragraph{Computation of torsion via disconnected spiders} 
In this paragraph we show that 
to determine the family of lifted cells necessary to define torsion one can use representatives
of Euler structures more general than those used above. This is a technical point which we
will use below to compute torsions using branched spines (Section~\ref{spines:section}).

We fix $M$, $\pp$, $\cc$ and $\varphi$ as above, and $\xi\in\eulc(M,\pp)$. Let
$\gbasis(\xi)=\{\tilde\sigma\}$ be the family of liftings of the cells lying in
$M\setminus(W\cup V)$ determined by a connected spider as explained above. Note
that if $\gbasis'=\{\tilde\sigma'\}$ is any other family of liftings we have 
$\tilde\sigma'=a(\sigma)\cdot\tilde\sigma$ for some $a\in\pi_1(M)$, and we can
define $$h(\gbasis',\gbasis(\xi))=\sum_{\sigma\subset M\setminus (W\cup
V)}\index(\sigma)\cdot \overline{a}(\sigma)\in H_1(M;\mz).$$

\begin{prop}\label{no:need:to:lift}
Assume there exists a partition $\cc_1\sqcup\dots\sqcup\cc_k$ of the set of cells 
lying in $M\setminus(W\cup V)$, and let $\xi\in\eulc(M,\pp)$ have a representative of the form
$$z=\sum_{j=1}^k\left(\sum_{\sigma\in\cc_j\setminus\{\sigma_j\}}
\index(\sigma)\cdot \gamma^{(j)}_\sigma\right)$$
where $\sigma_j\in\cc_j$ and $\gamma^{(j)}_\sigma:([0,1],0,1)\to(M,p_{\sigma_j},p_\sigma)$.
Choose any lifting $\tilde p_{\sigma_j}$ of $p_{\sigma_j}$, lift 
$\gamma^{(j)}_\sigma$ to $\tilde\gamma^{(j)}_\sigma$ starting from $\tilde p_{\sigma_j}$,
let $\tilde\sigma'$ be the lifting of $\sigma$ containing $\tilde\gamma^{(j)}_\sigma(1)$,
and let $\gbasis'$ be the family of all these liftings.
Then $h(\gbasis',\gbasis(\xi))=0\in H_1(M;\mz)$. In particular $\gbasis'$ can be used to compute 
$\tau^\varphi(M,\pp,\xi,\hbasis)$.
\end{prop}

\dim{no:need:to:lift}
Note first that the coefficient of $p_{\sigma_j}$ in $\partial z$ is exactly
$$-\sum_{\sigma\in\cc_j\setminus\{\sigma_j\}}\index(\sigma).$$
On the other hand this coefficient must be equal to $\index(\sigma_j)$. 
Summing up we deduce that $\sum_{\sigma\in\cc_j}\index(\sigma)=0$.

Now choose $x_0\in M$ and $\delta^{(j)}:([0,1],0,1)\to(M,x_0,p_{\sigma_j})$. For 
$\sigma\in\cc_j$ define
$$\beta_\sigma=\cases{
\delta^{(j)} & if $\sigma=\sigma_j$ \cr
\delta^{(j)}\cdot\gamma^{(j)}_\sigma & otherwise,}$$
so that $\beta_\sigma:([0,1],0,1)\to(M,x_0,p_\sigma)$, whence 
$w=\sum_{\sigma\subset M\setminus(W\cup V)}\beta_\sigma$ is an Euler chain.
Moreover:
$$w-z=\sum_{j=1}^k\left(\sum_{\sigma\in\cc_j}\index(\sigma)\right)
\cdot\delta^{(j)}=0\in H_1(M;\mz),$$
so $[w]=\xi$. Now choose $\tilde x_0$ over $x_0$, lift the $\delta^{(j)}$ and 
$\beta_\sigma$ starting from $\tilde x_0$, and let
$a^{(j)}\in\pi_1(M)$ be such that $\tilde p_{\sigma_j}=a^{(j)}\cdot\tilde\delta^{(j)}(1)$.
Then
$$h(\gbasis',\gbasis(\xi))=\sum_{j=1}^k\left(\sum_{\sigma\in\cc_j}\index(\sigma)\right)
\cdot\overline{a}^{(j)}=0\in H_1(M;\mz),$$
and the proof is complete.\finedim{no:need:to:lift}

The next result follows directly from the definition, but it is worth stating
 because it shows how torsions may be used to distinguish triples $(M,\pp,\xi)$ from each other.

\begin{prop}\label{homeo:action}
Let $f: M\to M'$ be a homeomorphism, consider $\xi\in\eul(M,\pp)$, 
$\varphi:\pi_1(M)\to\Lambda_*$ and a $\Lambda$-basis $\hbasis$ of 
$H^\varphi_*(M,\overline{W})$. Then
$$\tau^{\varphi\circ f_*^{-1}}(M',f_*(\pp),f_*(\xi),f_*(\hbasis))=
\tau^\varphi(M,\pp,\xi,\hbasis).$$
\end{prop}

\section{Spines and computation of torsion}\label{spines:section}
In this section we show how to geometrically invert the reconstruction map $\Psi$, and
how to compute torsions starting from a combinatorial
encoding of vector fields. We first review
the theory developed in~\cite{lnm}. See the beginning of Section~\ref{Eul:def:section}
for our conventions on manifolds, maps, and fields. In addition to the terminology
introduced there, we will need the notion of {\em traversing} field on a manifold
$M$, defined as a field whose orbits eventually intersect $\partial M$ 
transversely in both directions (in other words, orbits are compact intervals).

\paragraph{Branched spines}
A {\em simple} polyhedron $P$ is a  finite connected 2-dimensional polyhedron
with singularity of stable nature (triple lines and points where six non-singular
components meet). Such a $P$ is called {\it standard} if all the components of
the natural stratification given by singularity are open cells. Depending on
dimension, we will call the components {\it vertices, edges} and {\it regions}.

\begin{figure}
\centerline{\psfig{file=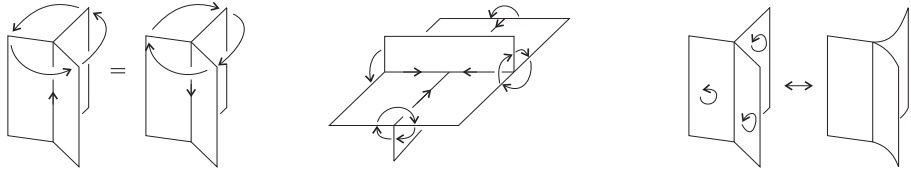,width=12cm}}
\caption{\label{screw:branch} Convention on screw-orientations, compatibility 
at vertices, and geometric interpretation of branching.}
\end{figure}

A {\em standard spine} of a $3$-manifold $M$ with $\partial M\neq\emptyset$
is a standard polyhedron $P$ embedded in $\interior(M)$ so that $M$ collapses onto $P$.
Standard spines of oriented $3$-manifolds are characterized among standard polyhedra
by the property of carrying an {\em orientation}, defined 
(see Definition~2.1.1 in~\cite{lnm}) as a ``screw-orientation''
along the edges (as in the left-hand-side of Fig.~\ref{screw:branch}),
with an obvious compatibility at vertices
(as in the centre of Fig.~\ref{screw:branch}).
It is the starting point of the theory of standard spines
that every oriented $3$-manifold $M$ with $\partial M\neq\emptyset$ has an oriented
standard spine, and can be reconstructed (uniquely up to homeomorphism) 
from any of its oriented standard spines. See~\cite{casler} for the non-oriented version 
of this result and~\cite{manuscripta} or Proposition~2.1.2 in~\cite{lnm} for the (slight) oriented refinement.

A {\it branching} on a standard polyhedron $P$ is an
orientation for each region of $P$, such that no edge is induced the same
orientation three times. See the right-hand side of Fig.~\ref{screw:branch}
and Definition 3.1.1 in~\cite{lnm} for the geometric meaning of this notion.
An oriented standard spine $P$ endowed with a branching is shortly named 
{\em branched spine}. We will never use specific notations for the extra structures: 
they will be considered to be part of $P$.
The following result, proved as Theorem~4.1.9 in~\cite{lnm}, is the starting point of our constructions.

\begin{prop}\label{from:spine:to:field}
To every branched spine $P$ there corresponds a manifold $M(P)$ 
with non-empty boundary and a concave traversing field $v(P)$ on $M(P)$.
The pair $(M(P),v(P))$ is well-defined up to diffeomorphism.
Moreover an embedding $i:P\to\interior(M(P))$ is defined, 
and has the property that $v(P)$ is positively transversal to $i(P)$.
\end{prop}

\begin{figure}
\centerline{\psfig{file=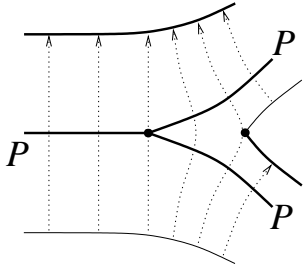,width=4cm}}
\caption{\label{constr:M} Manifold and field associated to a branched spine.}
\end{figure}

The topological construction which underlies this proposition is actually quite
simple, and it is illustrated in Fig.~\ref{constr:M}. Concerning the last
assertion of the proposition, note that the branching allows to define an
oriented tangent plane at each point of $P$.

\paragraph{Combinatorial encoding of combings}
Let $P$ be a branched spine, and define $v(P)$ on $M(P)$ as just explained.
Assume that in $\partial M(P)$ there is only one component which is homeomorphic to $S^2$
and is split by the tangency line of $v(P)$ to $\partial M(P)$ into two
discs. (Such a component will be denoted by $\stwotriv$.) Now, notice that $\stwotriv$ is
also the boundary of the closed $3$-ball with constant vertical field, 
denoted by $B^3_{\rm triv}$. This shows that we can cap off
$\stwotriv$ by attaching a copy of $B^3_{\rm triv}$, 
getting a compact manifold $\hatM(P)$ and a field $\hatv(P)$ on $\hatM(P)$. If we denote by $\hatpp(P)$ the boundary pattern of $\hatv(P)$ on
$\hatM(P)$, we easily see that the pair $(\hatM(P),\hatv(P))$ is only well-defined
up to homeomorphism of $\hatM(P)$ and homotopy of $\hatv(P)$ through fields compatible
with $\hatpp(P)$. Note also that $\hatpp(P)$ is automatically concave.

If $\pp$ is a boundary pattern on $M$, we define ${\rm Comb}(M,\pp)$ as
the set of fields compatible with $\pp$ under homotopy through fields also
compatible with $\pp$. An element of ${\rm Comb}(M,\pp)$
is called a {\em combing} on $(M,\pp)$. Note that we have a projection
${\rm Comb}(M,\pp)\to\eul(M,\pp)$.

The above construction shows that
a branched spine $P$ with only one $\stwotriv$ on $\partial M(P)$ defines
an element $\Phi(P)$ of ${\rm Comb}(\hatM(P),\hatpp(P))$.
In~\cite{second:paper} we will establish the following:

\begin{teo}\label{bounded:surg} 
If $M$ is any compact oriented 
$3$-manifold and $\pp$ is a concave boundary pattern
on $M$ not containing $\stwotriv$ components, then $\Phi$ maps surjectively
$\{P:\ \hatM(P)\cong M,\ \hatpp(P)\cong\pp\}$ onto ${\rm Comb}(M,\pp)$.
\end{teo}

This theorem generalizes
the main achievement of~\cite{lnm} (Theorems~1.4.1 and~5.2.1),
where it is proved in the special case of closed $M$. The complete statement
includes also the description of a finite set of local moves on branched spines
generating the equivalence relation induced by $\Phi$. We will not need the
moves in this paper. The following
geometric interpretation the theorem may however be of some interest.

\begin{rem}\label{trivial:pieces}
{\em In general, the dynamics of a field, even a concave one, can be very complicated,
whereas the dynamics of a traversing field (in particular, $B^3_{\rm triv}$) is simple.
Theorem~\ref{bounded:surg} means that for any 
(complicated) concave field
there exists a sphere $S^2$ which splits the field into two (simple) pieces:
a standard $B^3_{\rm triv}$ and a concave traversing field.}
\end{rem}

We can give here an easy special proof of Theorem~\ref{bounded:surg}
for the case we are most interested in, namely link exteriors.
Note that our argument relies on the results of~\cite{lnm}.

\vspace{1pt}
\noindent{\em Proof of~\ref{bounded:surg} 
for link exteriors}. We have to show that if $M$ is closed, $v$ is a field on $M$
and $L$ is transversal to $v$, then the exterior $E(L)$ of $L$ with the restricted field
is represented by some branched spine in the sense explained above.

The construction explained in Section~5.1  of~\cite{lnm} shows that 
there exists a branched standard spine $P$ such that $v$ is positively 
transversal to $P$, and the complement of $P$, with the 
restriction of $v$, is isomorphic to
the open 3-ball with the constant vertical field. 
The last condition easily implies that $L$ can be 
isotoped through links transversal to $v$ to a link lying in an arbitrarily small 
neighbourhood of $P$, with the further property that its natural projection
on $P$ is $\cont^1$, possibly with crossings. 

\begin{figure}
\centerline{\psfig{file=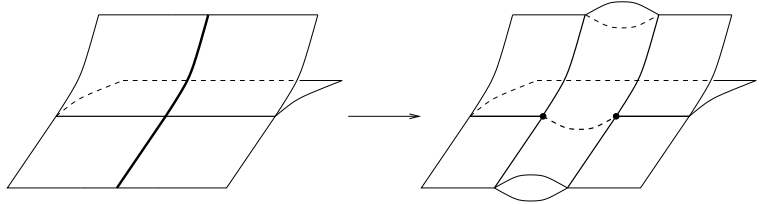,width=10cm}}
\caption{\label{dig:tunnel} How to dig a tunnel in a spine.}
\end{figure}

Once $L$ has been isotoped to a $\cont^1$ link on $P$,
a branched spine of $(E(L),v\ristr{E(L)})$ is obtained
by digging a tunnel in $P$ along the projection of $L$, as shown in Fig.~\ref{dig:tunnel}.
A crossing in the projection will of course give rise to 4 vertices in the spine.
Note that the spine which results from the digging may occasionally be 
non-standard, but it is standard as soon as the projection is complicated enough 
({\em e.g.}~if on each component there are both a crossing and an intersection with $S(P)$).
{{\hfill\hbox{\enspace\fbox{\ref{bounded:surg}$E(L)$}}}\vspace{5pt}}

\paragraph{Spines and ideal triangulations} 
We remind the reader that an {\it
ideal triangulation} of a manifold $M$ with non-empty boundary is  a partition
$\calt$ of $\interior(M)$ into open cells of dimensions 1, 2 and 3, induced by a
triangulation $\calt'$ of the space $Q(M)$, where:
\begin{enumerate}
\item $Q(M)$ is obtained from $M$ by collapsing each component of
$\partial M$ to a point;
\item $\calt'$ is a triangulation only in a loose sense, namely
self-adjacencies and multiple adjacencies of tetrahedra are allowed;
\item The vertices of $\calt'$ are precisely the points of $Q(M)$
which correspond to components of $\partial M$.
\end{enumerate}

\begin{figure}
\centerline{\psfig{file=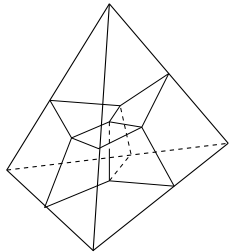,width=3cm}}
\caption{\label{duality} Duality between standard spines and ideal triangulations.}
\end{figure}

\noindent It turns out (see for instance~\cite{mafo},~\cite{tesi},~\cite{matv:new}) that there exists a natural
bijection between standard spines and ideal triangulations
of a 3-manifold. Given an ideal triangulation, the
corresponding standard spine is just the 2-skeleton of the dual
cellularization, as illustrated in Figure~\ref{duality}.
The inverse of this correspondence will be denoted by $P\mapsto\calt(P)$. 

Now let $P$ be a branched spine. First of all, we can realize $\calt(P)$ in such
a way that its edges are orbits of the restriction of $v(P)$ to
$\interior(M(P))$, and the 2-faces are unions of such orbits.  Being orbits,
the edges of $\calt(P)$ have a natural orientation, and the branching condition,
as remarked in~\cite{gr:1},
is equivalent to the fact that on each tetrahedron of $\calt(P)$ 
exactly one of the vertices is a sink and one is a source. 

\begin{rem}\label{all:oriented}
{\em It turns out that if $P$ is a branched spine, 
not only the edges, but also the faces
and the tetrahedra of $\calt(P)$ have natural orientations. For tetrahedra, we
just restrict the orientation of $M(P)$. For faces, we first note that the edges
of $P$ have a natural orientation (the prevailing orientation induced by the
incident regions). Now, we orient a face of $\calt(P)$ so that the algebraic
intersection in $M(P)$ with the dual edge is positive. }
\end{rem}

\paragraph{Euler chain defined by a branched spine}
We fix in this paragraph a standard spine $P$ and consider its manifold $M=M(P)$.
We start by noting that the ideal triangulation $\calt=\calt(P)$ defined by $P$ 
can be interpreted as a realization of $\interior(M)$ by face-pairings on a finite
set of tetrahedra with vertices removed. If, instead of 
removing vertices, we remove open conic neighbourhoods of the vertices, thus getting {\em truncated} tetrahedra, after the 
face-pairings we obtain $M$ itself. This shows that $P$ determines a cellularization
$\caltbar=\caltbar(P)$ 
of $M$ with vertices only on $\partial M$ and 2-faces which are either triangles
contained in $\partial M$ or hexagons contained in $\interior(M)$, with edges contained alternatingly in $\partial M$ and in $\interior(M)$.

\begin{figure}
\centerline{\psfig{file=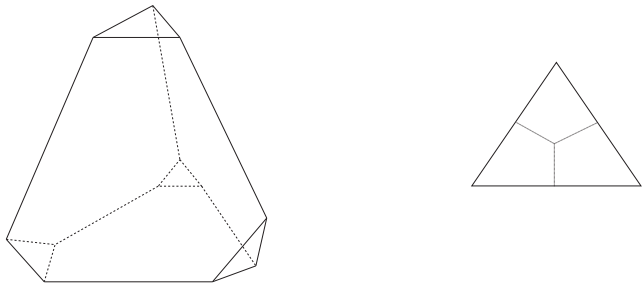,width=8.5cm}}
\caption{\label{trunc:tetra}Truncated tetrahedra and subdivision of the triangles on the
boundary}
\end{figure}

Now assume that $P$ is branched and 
that $\partial M$ contains only one $\stwotriv$ component, so 
$\hatM=\hatM(P)$ is defined. Note that $\hatM$ can be thought of as the space obtained from
$M$ by contracting $\stwotriv$ to a point, so a projection $\pi:M\to\hatM$ is defined,
and $\pi(\caltbar)$ is a cellularization of $\hatM$. Next, we modify 
$\pi(\caltbar)$ by subdividing the triangles on
$\partial\hatM$ as shown in Fig.~\ref{trunc:tetra}.
The result is a cellularization $\hattt=\hattt(P)$ of $\hatM$. 
Note that $\hattt$ on $\partial\hatM$ 
consists of ``kites'', with long edges coming from tetrahedra and short edges coming from 
subdivision. Note also that $\hattt$ has exactly one vertex $x_0$ in $\interior(\hatM)$, and
that the cells contained in $\interior(\hatM)$, except $x_0$, are the duals to the 
cells of the natural cellularization $\uu=\uu(P)$ of $P$. For $u\in\uu$ we denote by $\hat{u}$ 
its dual and by $p_u=p_{\hat{u}}$ the point where $u$ and $\hat{u}$ intersect, called the {\em centre} of both.

We will now use the field $\hatv=\hatv(P)$ to construct a combinatorial Euler chain
on $\hatM$ with respect to
$\hattt$. It is actually convenient to consider, instead of $\hatv$, the field 
$\vbar=\pi(v)$, which 
coincides with $\hatv$ except near $x_0$, where it has a 
(removable) singularity. For $u\in\uu$ we denote by $\beta_u$ the arc obtained by integrating $\vbar(P)$ in the positive direction, starting from $p_u$, until the 
boundary or the singularity is reached. We define:
$$s(P)=\sum_{u\in\uu}\index(u)\cdot \beta_u.$$

\begin{figure}
\centerline{\psfig{file=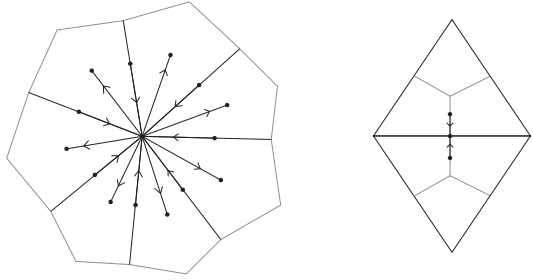,width=7cm}}
\caption{\label{starfig}The star ${\rm St}(p)$ centred at a vertex $p$ contained in $B$
and the bi-arrow ${\rm Ba}(\sigma)$ based at the midpoint of an edge $\sigma$ contained in $B$}
\end{figure}

Let us consider now the pattern $\hatpp=\hatpp(P)=(W,B,\emptyset,C)$ defined by $P$.
If $p$ is a vertex of $\pi(\caltbar)$ contained in $B$, we define its star ${\rm St}(p)$ as
the sum of the straight segments going from $p$ to the centres of all the kites containing $p$,
minus the sum of the straight
segments going from $p$ to the centres of all the long edges 
containing $p$. If $\sigma$ is an edge of $\pi(\caltbar)$ contained in $B$ we define its
bi-arrow ${\rm Ba}(\sigma)$ as the sum of the two straight segments going from the centre
$p_\sigma$ of $\sigma$ to the centres of the two short kite-edges containing $p_\sigma$.
A star and a bi-arrow are shown in Fig.~\ref{starfig}.
We define:
$$s'(P)=s(P)+\sum_{p\in B\cap\caltbar(P)^{(0)}}{\rm St}(p)+
\sum_{\sigma\in\hattt(P)^{(1)},\sigma\subset B}{\rm Ba}(\sigma).$$

\begin{lem}\label{s:prime:Euler}
$s'(P)$ defines an element of $\eulc(\hatM,\theta(P))$.
\end{lem}

\dim{s:prime:Euler}
Recall that $\theta(\hatpp)=(W,B,C,\emptyset)$, {\em i.e.}~the concave line $C$ is turned into
a convex one. So by definition we have to show that $\partial s'(P)$ contains, with the right 
sign, the centres of all cells of $\hattt$ except those of $W\cup C$. 

It will be convenient to analyze first the natural lifting of $s(P)$ to $M$, denoted by
$\tilde s(P)=\sum_{u\in\uu}\index(u)\cdot \tilde\beta_u$ with obvious meaning of symbols. So 
\begin{equation}
\partial\tilde s(P)=\sum_{u\in\uu}-\index(u)\cdot \tilde\beta_u(0)+
\sum_{u\in\uu}\index(u)\cdot \tilde\beta_u(1).
\label{s:tilde:boundary}
\end{equation}

Since the cellularization $\caltbar$ of $M$ is dual to $\uu$, the first half of 
(\ref{s:tilde:boundary}) gives the centres of the cells contained in $\interior(M)$, with right 
sign. One easily sees that the second half gives exactly the centres of the cells 
(of $\caltbar$) contained in $B$, also with right sign. 

When we project to $\hatM$ and consider $\partial s(P)$, the first half of 
(\ref{s:tilde:boundary}) again provides (with right sign)
the centres of the all cells contained in $\interior(\hatM)$, except the special
vertex $x_0$ obtained by collapsing $\stwotriv$. We can further split
the points of the second half of (\ref{s:tilde:boundary}) into those which lie on $\stwotriv$
and those which do not. The points of the first type project to $x_0$, and the resulting coefficient of $x_0$ is $\chi(B\cap\stwotriv)$, but $B\cap\stwotriv$ is an open 2-disc,
so the coefficient is 1. (We are here using the very special property of dimension 2
that $\chi$ can be computed using a finite cellularization of an open manifold, because the
boundary of the closure has $\chi=0$.) The points of the second type faithfully project to 
$\hatM$, giving the centres of the simplices contained in $B$ of the triangulation
$\pi(\caltbar)\ristr{\partial\hatM}$.
However $\hattt$ on $\partial\hatM$ is a subdivision of $\pi(\caltbar)$,
and this is the reason why we have added the stars and the bi-arrows to $s(P)$ getting $s'(P)$.
The following computation of the coefficients in $\partial s'(P)$ of the centres of the 
cells of $\hattt$ contained in $B$ concludes the proof.
\begin{enumerate}
\addtocounter{enumi}{-1}
\item Cells of dimension 0 are listed as follows:
\begin{enumerate}
\item Centres of triangles of $\pi(\caltbar)$, 
which receive coefficient $+1$ from $\partial s(P)$;
\item Midpoints of edges of $\pi(\caltbar)$, 
which receive coefficient $-1$ from $\partial s(P)$
and $+2$ from the bi-arrows they determine;
\item Vertices of $\pi(\caltbar)$ receive $+1$ from $\partial s(P)$ and
(algebraically) $0$  from the star they determine;
\end{enumerate}
\item Cells of dimension 1 are:
\begin{enumerate}
\item Short edges of kites, whose midpoints receive $-1$ from the bi-arrows;
\item Long edges of kites, whose midpoints receive $-1$ from the stars;
\end{enumerate}
\item Cells of dimension 2 are kites, and their centres receive $+1$ from the stars.
\end{enumerate}
\finedim{s:prime:Euler}

\noindent Now we denote by $\gamma_j:(0,1)\to C$, for $j=1,\dots,n$, 
orientation-preserving parameterizations of the 1-cells of $\hattt$ contained in $C$, and we extend the $\gamma_j$ to $[0,1]$, without changing notation. We define
$$s''(P)=s'(P)+\sum_{j=1}^n\gamma_j\ristr{[1/2,1]}.$$ 

\begin{lem}\label{s:second:Euler}
$s''(P)$ defines an element of $\eulc(\hatM,\hatpp)$, and 
$$[s'(P)]=\Thetac([s''(P)])\in \eulc(\hatM,\theta(\pp)).$$
\end{lem}

\dim{s:second:Euler} At the level of representatives, the second assertion is obvious,
and it implies the first assertion.
\finedim{s:second:Euler} 

\noindent We defer to Section~\ref{proofs} the proof of the next result, which shows that
the map $P\mapsto[s''(P)]\in\eulc(\hatM,\hatpp)$ allows,
using branched spines, to explicitly find the inverse
of the reconstruction map $\Psi$ of 
Theorem~\ref{reconstruction:statement}. This result was informally announced as
Theorem~\ref{informal:algorithmic} in the Introduction.

\begin{teo}\label{spider:structure:teo}
$\Psi([s''(P)])=[\hatv(P)]\in\euls(\hatM,\hatpp)$.
\end{teo}

Recall now that we have defined torsions directly only for convex patterns,
and we have extended the definition to concave patterns via the map $\Theta$.
As a consequence of Lemma~\ref{s:second:Euler} and Theorem~\ref{spider:structure:teo},
and by direct inspection of $s'(P)$, we have the following result which summarizes
our investigations on the relation between spines, Euler structures, and torsion:

\begin{teo}\label{spines:compute:teo}
If $P$ is a branched spine which represents a manifold $\hatM$ with 
concave boundary pattern $\hatpp=(W,B,\emptyset,C)$
in the sense of Theorem~\ref{bounded:surg}, then for any 
representation $\varphi:\pi_1(M)\to\Lambda_*$ and any $\Lambda$-basis
$\hbasis$ of $H^\varphi_*(\hatM,W\cup C)$, the torsion
$\tau^\varphi(\hatM,\hatpp,[\hatv(M)],\hbasis)$ can be computed using
(in the sense of Proposition~\ref{no:need:to:lift}) the lifting
to the universal cover of $\hatM$ of the chain $s'(P)$ defined above. In particular, 
$s'(P)$ can be used directly, without replacing it by a connected spider.
\end{teo}

\paragraph{Computational hints} 
To actually compute torsion starting from a branched spine $P$, besides describing
the universal (or maximal Abelian) cover of $\hatM=\hatM(P)$ and determining the preferred
liftings of the cells in $\hatM\setminus(W\cup C)$, one needs to compute
the boundary operators in the twisted chain complex 
$\Cphi_*(M,W\cup C)$. These operators are of course twisted liftings of the
corresponding operators in the cellular chain complex of $(\hatM,W\cup C)$,
with respect to $\hattt$. We briefly describe here the form of the latter operators.
Recall first that $\hattt$ consists of a special vertex $x_0$, the kites
(with their vertices and edges) on $\hatM$, and the duals of the
cells of $P$. On $\partial\hatM$ the situation is easily described, so we consider the
internal cells. 
\begin{enumerate}
\item If $R$ is a region of $P$, the ends of its dual edge 
$\hat{R}$ are either $x_0$ or vertices of $\partial\hatM$ 
contained only in long edges of kites. 
\item If $e$ is an edge of $P$ then
$\partial\hat{e}$ is given by
$\hat{R}_1+\hat{R}_2-\hat{R}_0$ plus 3 long edges of kites, where $R_0,R_1,R_2$ are the
regions incident to $e$, numbered so that $R_1$ and $R_2$ induce on $e$ the same
orientation. Here $R_0,R_1,R_2$ need not be different from each other, so
the formula may actually have some cancellation. The 3 long edges
of kites must be given an appropriate sign, and some of them 
may actually be collapsed to the point $x_0$. Note that we have only 3 kite-edges, out
of the 6 which geometrically appear on $\partial\hat{e}$, because the other 3 are white.
\item If $v$ is a vertex of $P$ then
$\partial\hat{v}$ is given by 
$\hat{e}_1+\hat{e}_2-\hat{e}_3-\hat{e}_4$ plus 6 kites, where $e_1,e_2$ are
the edges which (with respect to the natural orientation) are leaving $v$, and
$e_3,e_4$ are those which are reaching it. Again, there could be
repetitions in the $e_i$'s. The kites all have coefficient $+1$, and again some of them may
actually be collapsed to $x_0$. As above, we have only
6 kites because the other 6 are white.
\end{enumerate}

\begin{rem}\label{do:not:cut}{\em To define the cellularization $\hattt(P)$ associated to 
a spine we have decided to subdivide all the triangles on $\partial\hatM$
into 3 kites, but when doing actual computations this is not necessary and 
impractical. The only triangles 
which we really need to subdivide are those intersected by $C$, because we need the 
cellularization to be suited to the pattern. Let us consider the 4
triangles corresponding to the ends of a certain tetrahedron. If in each of them
we count the number of black kites and the number of white kites, we get respectively
$(3,0)$, $(2,1)$, $(1,2)$, $(0,3)$. So, the first and last triangles do not have to
be subdivided, and the other two can be subdivided using a segment only. Summing up,
for each vertex of $P$ we only need to add two segments on the boundary. Before
projecting $M(P)$ to $\hatM(P)$ one sees that the number of cells, with respect to
$\caltbar(P)$, is increased in all dimensions 0, 1 and 2 by twice the number of vertices of $P$.
When projecting to $\hatM(P)$ the cells lying in $\stwotriv$ get collapsed to points.}
\end{rem}

\section{Pseudo-Legendrian knots}\label{knots:section}

We fix in this section a compact oriented manifold $M$ and a boundary pattern $\pp$ on $M$.
The boundary of $M$ may be empty or not. Recall that if $v$ is a vector field on $M$ and 
$K$ is a knot in $\interior(M)$, we have defined $K$ to be pseudo-Legendrian in $(M,v)$ 
if $v$ is transversal to $K$. We will also call $(v,K)$ a pseudo-Legendrian pair.
Having fixed $\pp$, we will only consider fields $v$ 
compatible with $\pp$. Some of the results we will establish hold
also for links, but we will stick to knots for the sake of 
simplicity. First, we need to spell out the equivalence relation
which we consider.

Let $v_0,v_1$ be compatible with $\pp$ and let $K_0,K_1$ be pseudo-Legendrian in
$(M,v_0)$ and $(M,v_1)$ respectively. We define $(v_0,K_0)$ to be {\em weakly equivalent}
to $(v_1,K_1)$ if there exist a homotopy $(v_t)_{t\in[0,1]}$ through fields compatible
with $\pp$ and an isotopy $(K_t)_{t\in[0,1]}$ such that $K_t$ is transversal
to $v_t$ for all $t$. If $v_0=v_1$ then $K_0$ and $K_1$ are
called {\em strongly equivalent} if the homotopy $(v_t)$ can be chosen to be constant.

\begin{rem}\label{weak:strong:rem}
{\em
\begin{enumerate}
\item Of course strong equivalence implies weak equivalence. Weak equivalence is the natural
relation to consider on pseudo-Legendrian pairs $(v,K)$, while strong equivalence
is natural for pseudo-Legendrian knots in a fixed $(M,v)$. See~\cite{second:paper}
for a further discussion on this point.
\item The term `pseudo-Legendrian isotopy', used in the introduction and 
in~\cite{second:paper}, corresponds to `weak equivalence.' For the sake of brevity,
and to emphasize the difference with strong equivalence, we will only use 
the term `weak equivalence' in the rest of this paper.
\end{enumerate}}
\end{rem}

Before proceeding note that if $K$ is pseudo-Legendrian in $(M,v)$ then
$v$ turns $K$ into a framed knot, which we will denote by $K^{(v)}$.
The framed-isotopy class of $K^{(v)}$ is of course invariant under weak
equivalence.

\paragraph{Euler structures on knot exteriors}
For a knot $K$ in $M$ we consider a (closed) tubular neighbourhood $U(K)$ of $K$
in $M$ and we define $E(K)$ as the closure of the complement of $U(K)$.
If $F$ is a framing on $K$ we extend the boundary pattern $\pp$ previously fixed
on $M$ to a boundary pattern $\pp(K^F)$ on $E(K)$, by splitting $\partial U(K)$
into a white and a black longitudinal annuli, the longitude being the one defined by
the framing $F$. As a direct application of Proposition~\ref{p:h:formula} one sees
that $\eul(E(K),\pp(K^F))$ is non-empty (assuming $\eul(M,\pp)$ to be non-empty).

A convenient way to think of $\pp(K^F)$ is as follows. The framing $F$ 
determines a transversal vector field along $K$. If we extend this field near
$K$ and choose  $U(K)$ small enough then the pattern we see on $\partial U(K)$
is exactly as required. With this picture in mind, it is clear that if $K$ is
pseudo-Legendrian in $(M,v)$, where $v$ is compatible with $\pp$, then the
restriction of $v$ to $E(K)$ defines an element
$$\xi(v,K)\in\eul(E(K),\pp(K^{(v)})$$ 
(this structure was already considered in the partial proof of Theorem~\ref{bounded:surg}).
So the theory of torsion applies. We will discuss in this section torsion both
as an {\em absolute} and as a {\em relative} invariant, splitting the section into
two subsections.

\subsection{Absolute torsion and the Alexander invariant}
We will show in this section that torsion as an absolute invariant of
pseudo-Legendrian knots lifts the classical Alexander invariant, but the
relation between the two objects is more complicated than in Turaev's
situation (\cite{turaev:Reidemeister} and~\cite{turaev:Euler}),
because here two different algebraic complexes are involved at the same time.

\paragraph{Turaev's lifting of Milnor torsion} 
Let us first recall again in what sense Turaev's torsion lifts the classical one. Let
$M$ be a manifold which is closed or bounded by tori, and take
a representation $\varphi:H_1(M;\mz)\to\Lambda_*$, where $\Lambda$
is as usual. The classical theory~\cite{milnor} allows to
define an invariant
$$\tau^\varphi(M)\in K_1(\Lambda)
\Big/ (\pm\varphi(H_1(M;\mz))),$$
usually stipulated to be $1$ if
the $\varphi$-twisted homology of $M$ does not vanish, {\em i.e.},
using the above notation, if the complex 
$C^\varphi_*(M)=\Lambda\otimes_\varphi
C_*^{\rm cell}(\tilde M;\mz)$ is not acyclic, where
$q:\tilde M\to M$ is the maximal Abelian cover.
When $\xi$ is an Euler structure on $M$ with monochromatic boundary components,
Turaev~\cite{turaev:Euler} 
shows that his torsion $\tau^\varphi(M,\xi)\in\Kbar_1(\Lambda)$ is a lifting
of $\tau^\varphi(M)$ with respect to the obvious projection of
$K_1(\Lambda)/\pm\varphi(H_1(M;\mz)$ onto $\Kbar_1(\Lambda)$.

In the special case where 
$\Lambda$ is the field of fractions obtained from the 
group ring of $H_1(M;\mz)$ modulo torsion,
and $\varphi:H_1(M;\mz)\to\Lambda$ is the natural projection, the invariant
$\tau^\varphi(M)$ is called Milnor torsion, and its
sign-refinement provided by Turaev 
in~\cite{turaev:Reidemeister} has been shown to be equivalent
to the classical Alexander invariant. So Turaev's torsion for Euler structures contains
a lifting of the Alexander invariant. We will discuss in the rest of this
subsection the extent to what the same holds when the Euler structure 
arises from a pseudo-Legendrian knot. What we will say applies to any allowed
representation $\varphi:H_1(M;\mz)\to\Lambda$, but we keep in mind that the relation
with the Alexander invariant emerges for a special choice of $\varphi$ and $\Lambda$.
Since we will also drop the condition that the involved complexes be acyclic, we note
that torsion is only defined when the resulting homology is free. This is not true in
general, but it is for instance when $\Lambda$ is a field.

\paragraph{Torsion of a knot complement} Let us restrict to the case of a closed
manifold $M$, and let us consider a pseudo-Legendrian pair $(v,K)$ in $M$
and a representation $\varphi:H_1(E(K);\mz)\to\Lambda$ as
usual. We would like to interpret the torsion of the Euler structure $\xi(v,K)$ 
on $E(K)$ with respect to $\varphi$ 
as a lifting of $\tau^\varphi(E(K))$, but a difficulty immediately emerges, because
the algebraic complexes used to compute these torsions do not coincide.

To be more specific, let us first spell out how the torsion of $\xi(v,K)$ is defined.
Let $\pp(K^{(v)})=(B,W,\emptyset,C)$ be the boundary pattern
defined on $E(K)$. Then we define
$\tau^\varphi(M,v,K,\hbasis)$ as $\tau^\varphi(M,\pp(K^{(v)}),\Theta(\xi(v,K)),\hbasis)$.
More specifically, 
$\tau^\varphi(M,v,K,\hbasis)$ is the torsion of the complex
$C^\varphi_*(E(K),\overline{W})$,
where $W$ is the (open) white annulus on $\partial E(K)$,
as above the maximal Abelian cover of $E(K)$ is used to define the complex, 
the preferred cell lifting is obtained using
an Euler chain for the convexified structure $\Theta(\xi(v,K))$, and $\hbasis$ is
a basis of the twisted homology of $E(K)$ relative to $\overline{W}$.

Now, $\tau^\varphi(E(K))$ is the torsion of
$C_*^\varphi(E(K))$, and
this complex can be radically different from the previous one. For instance,
when $M$ is a homology sphere, the absolute complex is always acyclic, while
the complex relative to $\overline{W}$, which depends only on
the framed knot $K^{(v)}$, very often is not. 
We will see how to overcome this difficulty
using the fundamental multiplicativity properties of torsion.

\begin{figure}
\centerline{\psfig{file=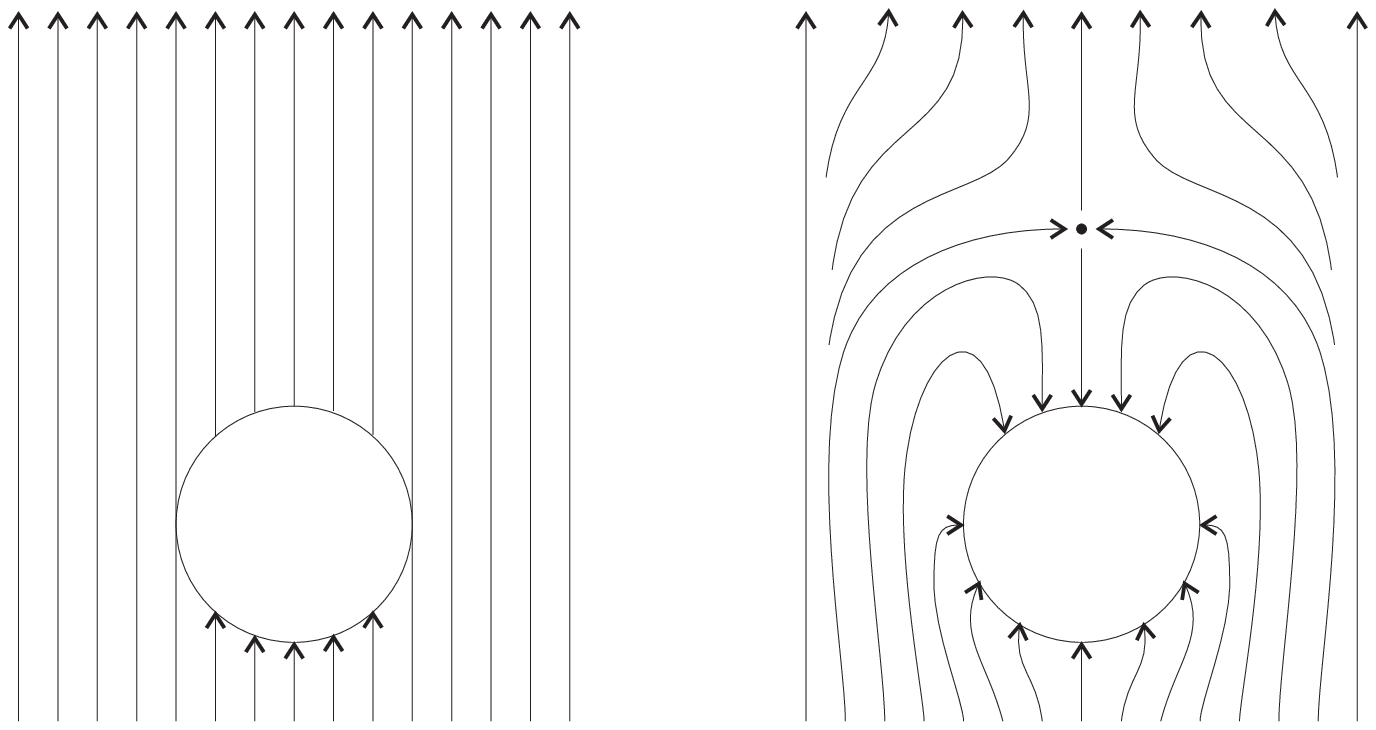,width=8cm}}
\caption{\label{geomblak} Black field on a knot complement.}
\end{figure}

\paragraph{How to turn a torus into black}
We will describe in this paragraph two explicit methods for
modifying $\xi(v,K)$ to an Euler structure $\beta(v,K)$
such that $\partial E(K)$ becomes monochromatic black.
These methods are respectively a geometric and an algebraic one. 
The fact that they actually lead to the same result is
true but not very important, so we will omit the proof.
Both methods involve the choice of an orientation of $K$.
The first method is explained in a cross-section in Fig.~\ref{geomblak}
The cross-section is transversal to $K$, and the apparent singularity
of the modified field is removed by summing a field parallel to $K$ and
supported near the singularity ({\em cf.}~Fig.~\ref{conc:to:conv}, where a similar
method was used).

\begin{figure}
\centerline{\psfig{file=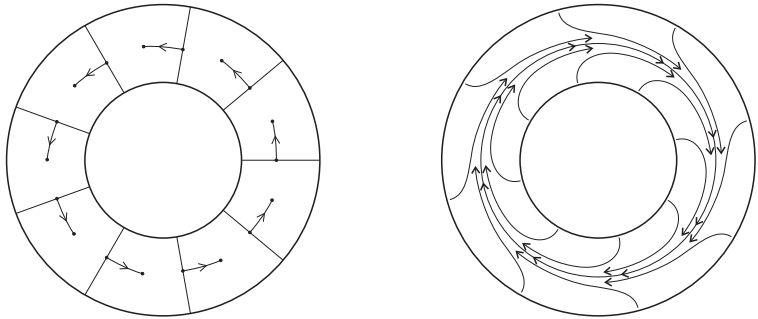,width=10cm}}
\caption{\label{zw_chain} A 1-chain on the annulus $W$.}
\end{figure}

To describe the algebraic construction of $\beta$, recall that if $z$ is a 1-chain
representing $\xi(v,K)$ then $\partial z$ contains,
with the appropriate sign, the centres of all cells in $E(K)\setminus W$. 
Knowing the subdivision rule for Euler chains (Proposition~\ref{combin:prop:statement})
we can also assume that the cellularization on $W$ has a particularly simple shape.
We assume it consists of rectangles as in Fig.~\ref{zw_chain} (left), where we also
show a 1-chain $z_W$ having the property that $\partial z_W$ contains the centres of all
cells in $W$. We can now define $\beta(v,K)$ as the Euler structure carried by $z+z_W$.
The boundary of $E(K)$ is completely black with respect to this structure,
because $\partial(z+z_W)$ contains the centres of all cells of $E(K)$.

One easily sees from both our descriptions of $\beta$ that it is canonically
defined and $H_1$-equivariant.
Since we will need these properties, we spell out their meaning, starting from 
an oriented framed knot $K^F$ rather than a pseudo-Legendrian knot.
Let $(W,B,\emptyset,C)$ be the concave boundary pattern determined
by $F$ on $E(K)$: then $\beta:\eul(E(K),(W,B,\emptyset,C))\to\eul(E(K);
(\emptyset,\partial E(K),\emptyset,\emptyset)$ 
is well-defined (depending on $K^F$ only) and $H_1(E(K);\mz)$-equivariant.

\begin{rem}{\em If $-K$ denotes the same knot $K$ with reversed orientation then
$$\alpha(\beta(v,K),\beta(v,-K))=[\lambda]\in H_1(E(K);\mz)$$
where $\lambda$ is the longitude on $\partial E(K)$ determined by the framing $K^{(v)}$.}
\end{rem}

A geometric interpretation of the second description of $\beta$ is possible
and used below. We have mentioned that a theory of Euler structures exists in all dimensions.
While the case $n\geq 4$ requires some technicalities~\cite{fourth:paper},
the reader can easily work out the case $n=2$ using the case
$n=3$ treated in the present paper. 
And one easily sees that $z_W$ is just an Euler chain of the inward-pointing
Reeb field $r_W$ on $W$ shown in Fig.~\ref{zw_chain} (right). Moreover $r_W$ can
be canonically turned into an outward-pointing field $\Theta(r_W)$, which of course
is the outward-pointing Reeb field (but the core spins in the opposite direction).
So a torsion $\tau^\psi(W,\Theta(r_W))$ can be computed
(possibly with a basis of the twisted homology added to the data).

\paragraph{Knot torsion as a lifting of Milnor's torsion}
Let as above $(v,K)$ be pseudo-Legendrian and let $\varphi:H_1(E(K);\mz)\to\Lambda$
be a representation. If $i:W\to E(K)$ is the inclusion, we set 
$\varphi_W=\varphi\compo i_*$. Considering the twisted homology of the pair
$(E(K),W)$ we get an exact sequence
$$\hh=\Big(\cdots\longrightarrow H^{\varphi_W}_i(W)\longrightarrow
H^\varphi_i(E(K))\longrightarrow H^\varphi(E(K),W)\longrightarrow
H^{\varphi_W}_{i-1}(W)\longrightarrow\cdots\Big).$$
We choose bases $\hbasis$, $\hbasis'$, and $\hbasis''$ respectively for 
$H^\varphi_*(E(K);\overline{W})$, $H^\varphi_*(E(K))$ and $H^{\varphi_W}_*(W)$,
so we can compute  
$\tau^\varphi(M,v,K,\hbasis)$, $\tau^\varphi(E(K),\beta(v,K),\hbasis')$ and 
$\tau^{\varphi_W}(W,\Theta(r_W),\hbasis'')$. Moreover we can compute 
$\tau(\hh,\hbasis,\hbasis',\hbasis'')$. The following
result is a refinement of Theorem~3.2 in~\cite{milnor}, and
a proof can be given imitating the argument given in~\cite{turaev:spinc}
(where a special case of the result is established).

\begin{prop}\label{Alex:lift}
The following equality holds:
\begin{equation}\label{Alex:lift:eq}
\tau^\varphi(E(K),\beta(v,K),\hbasis')=
\tau^\varphi(M,v,K,\hbasis)\cdot
\tau^{\varphi_W}(W,\Theta(r_W),\hbasis'')\cdot
\tau(\hh,\hbasis,\hbasis',\hbasis'').
\end{equation}
\end{prop}

The following comments on the previous proposition eventually explain
in what sense our torsion can be viewed as a lifting of the classical torsion
(in particular, Milnor torsion and the Alexander invariant).

\begin{rem}{\em 
In equation~(\ref{Alex:lift:eq}) the term $\tau^\varphi(E(K),\beta(v,K),\hbasis')$
is one of Turaev's torsion, so it is indeed a lifting of the classical torsion.
The term $\tau^\varphi(M,v,K,\hbasis)$ is the torsion for pseudo-Legendrian knots
introduced in this paper, while 
$\tau^{\varphi_W}(W,\Theta(r_W),\hbasis'')$ and 
$\tau(\hh,\hbasis,\hbasis',\hbasis'')$ can be viewed as normalizing terms.
One can for instance choose homology bases so that $\tau(\hh,\hbasis,\hbasis',\hbasis'')=1$,
and note that $\tau^{\varphi_W}(W,\Theta(r_W),\hbasis'')$ depends only
on the framed knot $K^{(v)}$, not on the Euler structure.}
\end{rem}

\begin{rem}{\em 
The factor $\tau^{\varphi_W}(W,\Theta(r_W),\hbasis'')$ 
may be understood quite easily. Denoting by $1$ the generator of $H_1(W;\mz)$,
the result depends on $\varphi_W(1)$. If
$\varphi_W(1)=1$ then the $\varphi_W$-twisted homology of $W$ is not
twisted at all, so it is non-zero and free, and we can
choose $\hbasis''$ so that $\tau^{\varphi_W}(W,\Theta(r_W),\hbasis'')=1$.
On the contrary, if $\varphi_W(1)-1$ is invertible, then the $\varphi_W$-twisted homology
is zero, and $\tau^{\varphi_W}(W,\Theta(r_W))$  is computed to be
$\varphi_W(1)-1$. In the intermediate cases where $\varphi_W(1)-1$ is neither zero nor a unit,
which can only occur when $\Lambda$ is not a field,
$\tau^{\varphi_W}(W,\Theta(r_W))$ is likely not to be defined.}
\end{rem}

\begin{rem}{\em
We can further specialize the understanding of $\tau^{\varphi_W}(W,\Theta(r_W),\hbasis'')$ 
when $M$ is a homology sphere and $\varphi:H_1(E(K);\mz)\to\Lambda$ is the representation
which gives the Milnor torsion. In this case $\varphi_W(1)=\varphi(\mu)^n$, where $\mu$ is the
meridian of $K$ and $n$ is the framing. So $\tau^{\varphi_W}(W,\Theta(r_W),\hbasis'')$
can be normalized to be $1$ for $\varphi(\mu)^n=1$, and it equals $\varphi(\mu)^n-1$ when
$\varphi(\mu)^n\neq1$}\end{rem}

\subsection{Torsion as a relative invariant}
We study in this section how torsion can be employed to distinguish
pseudo-Legendrian knots from each other. We first show that as a relative invariant
torsion is only well-defined as a multi-valued function, the ambiguity being
given by the action of a group. Then we concentrate on the knots (called `good' below)
for which this action is trivial, and we interpret the relative information carried
by torsion as a difference of winding numbers (or Maslov indices).

\paragraph{Group action on Euler structures}
Consider a knot $K$ and a
self-diffeomorphism $f$ of $E(K)$ which is the identity near $\partial E(K)$.
Then $f$ extends to a self-diffeomorphism $\widehat f$ of $M$, where 
$\widehat f\ristr{U(K)}={\rm id}_{U(K)}$. We define $G(K)$ as the group of all such $f$'s
with the property that $\widehat f$ is isotopic to the identity on $M$.
Elements of $G(K)$ are regarded up to isotopy relative to $\partial E(K)$.
If $F$ is a framing on $K$ then the pull-forward of vector fields induces
an action of $G(K)$ on $\eul(E(K),\pp(K^{(v)})$. We will now see 
that an obstruction to weak equivalence can be expressed in terms
this group action.

Let $(v_0,K_0)$ and $(v_1,K_1)$ be pseudo-Legendrian pairs in $M$, and assume
that  $K_0^{(v_0)}$ is framed-isotopic to $K_1^{(v_1)}$ under a diffeomorphism
$f$ relative to $\partial M$. Using the restriction of $f$ and the pull-back of
vector fields we get a bijection
$$f^*:\eul(E(K_1),\pp(K_1^{(v_1)}))\to\eul(E(K_0),\pp(K_0^{(v_0)})).$$

\begin{prop}\label{group:obstr}
Under the current assumptions, if $(v_0,K_0)$ and $(v_1,K_1)$ are weakly equivalent to each
other then $f^*(\xi(v_1,K_1))$ belongs to the $G(K_0)$-orbit of $\xi(v_0,K_0)$ in
$\eul(E(K_0),\pp(K_0^{(v_0)})$.
\end{prop}

\dim{group:obstr}
By assumption $K_0,K_1$ and $v_0,v_1$ embed in continuous families $(K_t)_{t\in[0,1]}$
and $(v_t)_{t\in[0,1]}$, where $v_t$ is transversal to $K_t$ for all $t$. Now
$(K_t^{(v_t)})_{t\in[0,1]}$ is a framed-isotopy, so 
there exists a continuous family $(g_t)_{t\in[0,1]}$ of diffeomorphisms of $M$
fixed on $\partial M$ and such that $g_0={\rm id}_M$ and $g_t(K_0^{(v_0)})=K_t^{(v_t)}$. 
So we get a map
$$[0,1]\ni t\mapsto
\alpha(\xi(v_0,K_0),g_t^*(\xi(v_t,K_t)))\in H_1(E(K_0);\mz).$$
Since $H_1(E(K_0);\mz)$ is discrete and the map is continuous,
we deduce that the map is identically 0. So $g_1^*(\xi(v_1,K_1))=\xi(v_0,K_0)$.
Now
$$f^*(\xi(v_1,K_1))=(f^*\compo (g_1)_*\compo g_1^*)(\xi(v_1,K_1))=
(f^{-1}\compo g_1)_*(\xi(v_0,K_0))$$
and the conclusion follows because $f^{-1}\compo g_1$ defines an element of $G(K_0)$.
\finedim{group:obstr}

The group $G(K)$ is in general rather difficult to understand (see~\cite{hatcher}), 
so we introduce a special terminology for the case where its action
can be neglected. We will say that a framed knot $K^F$ is {\em good} if
$G(K)$ acts trivially on $\eul(E(K),\pp(K^F))$. If $K^F$ is good for all framings $F$,
we will say that $K$ itself is good.
The following are easy examples of good knots:
\begin{itemize} 
\item $M$ is $S^3$ and $K$ is the trivial knot;
\item $M$ is a lens space $L(p,q)$ and $K$ is the core of one of the handlebodies
of a genus-one Heegaard splitting of $M$.
\end{itemize} 
\noindent The reason is that in both cases $E(K)$ is a solid torus, and we know that an automorphism of
the solid torus which is the identity on the boundary is isotopic 
to the identity relatively to the boundary, so $G(K)$ is trivial. The next three
results show that on one hand $G(K)$ is very seldom trivial, but on the other hand 
many knots are good. We will give proofs in the sequel, after introducing
some extra notation. In the statements, by `$E(K)$ is hyperbolic' 
we mean `$\interior(E(K))$ is complete, finite-volume hyperbolic.'

\begin{prop}\label{seldom:triv}
If $M$ is closed and $E(K)$ is hyperbolic then $G(K)$ is non-trivial.
\end{prop}

\begin{teo}\label{many:good}
If $M$ is closed, $E(K)$ is hyperbolic and 
either ${\rm Out}(\pi_1(E(K)))$ is trivial or $H_1(E(K);\mz)$ is torsion-free
then $K$ is good.
\end{teo}

\begin{teo}\label{homology:good}
If $M$ is a homology sphere then every knot in $M$ is good.
\end{teo}

The next result, which follows directly from Proposition~\ref{group:obstr},
the definition of goodness, and Proposition~\ref{homeo:action},
shows that for good knots torsion can be used as an
obstruction to weak (and hence strong) equivalence.

\begin{prop}\label{torsion:obstr}
Let $(v_0,K_0)$ and $(v_1,K_1)$ be pseudo-Legendrian pairs in $M$, and assume that 
$K_0^{(v_0)}$ is framed-isotopic to $K_1^{(v_1)}$ under a diffeomorphism $f$ relative
to $\partial M$. Suppose that $K_0^{(v_0)}$ is good, and that for some representation
$\varphi:\pi_1(E(K_0))\to\Lambda$ and some $\Lambda$-basis $\hbasis$ of
$H^\varphi_*(E(K_0),\overline{W(\pp(K_0^{(v_0)}))})$ we have
\begin{equation}\label{tors:obs:eq}
\tau^\varphi(E(K_0),\pp(K_0^{(v_0)}),\xi(v_0,K_0),\hbasis)\neq
\tau^{\varphi\compo f_*^{-1}}(E(K_1),\pp(K_1^{(v_1)}),\xi(v_1,K_1),f_*(\hbasis)).
\end{equation}
Then $(v_0,K_0)$ and $(v_1,K_1)$ are not weakly equivalent.
\end{prop}

\begin{rem}{\em
\begin{enumerate}
\item The right-hand side of equation~(\ref{tors:obs:eq}) actually equals
\begin{eqnarray*}
& & \tau^\varphi(E(K_0),\pp(K_0^{(v_0)}),f^*(\xi(v_1,K_1)),\hbasis) \\
&=& \overline{\varphi}(\alpha(v_0,f^*(v_1))\cdot
\tau^\varphi(E(K_0),\pp(K_0^{(v_0)}),\xi(v_0,K_0),\hbasis).
\end{eqnarray*}
This shows that the most torsion can capture as a relative invariant
of $(v_0,K_0)$ and $(v_1,K_1)$ is $\alpha(v_0,f^*(v_1))$. 
We will show below that in some cases torsion indeed allows
to determine $\alpha(v_0,f^*(v_1)$ completely.
\item By definition of goodness the homology class
$\alpha(v_0,f^*(v_1))$ just considered is actually
independent of $f$. We will denote it by $\alpha((K_0,v_0),(K_1,v_1))$.
\item For non-good knots the relative invariant is an orbit of the action of
$G(K_0)$. So an obstruction in terms of torsion could be given also for non-good knots,
but the statement would become awkward, and we have 
refrained from giving it.
\item If equation~(\ref{tors:obs:eq}) holds for some basis $\hbasis$ then it holds
for any basis.
\end{enumerate}}
\end{rem}

\noindent To conclude this paragraph we note that using the technology of
Turaev~\cite{turaev:Euler},  one can actually see that the action on Euler
structures of an automorphism  is invariant under {\em homotopy} (not only
isotopy) relative to the boundary. We will not use this fact.

\paragraph{Good knots} 
We introduce now some notation needed for the
proofs of Proposition~\ref{seldom:triv} and Theorem~\ref{many:good} (for 
Theorem~\ref{homology:good} we will use a different approach, see below).
Recall that $(M,\pp)$ is fixed for the whole section. We temporarily fix also
a framed knot $K^F$ in $M$, a regular neighbourhood $U$ of $K$, and we denote by $T$ the
boundary torus of $U$. On $T$ we consider 1-periodic coordinates $(x,y)$
such that $x\mapsto (x,0)$ is a meridian of $U$ and $y\mapsto (0,y)$ is 
a longitude compatible with $F$. We denote a collar of $T$ in $E(K)$ by 
$V$ and parameterize $V$ as $T\times[0,1]$, where $T=T\times\{0\}$. We consider on
$[0,1]$ a coordinate $s$. For $p,q\in\mz$ we define automorphisms
$\dd_{(p,q)}$ of $E(K)$ as follows. 
Each $\dd_{(p,q)}$ is supported in $V$, and on $V$, using the
coordinates just described, it is given by 
$$\dd_{(p,q)}(x,y,s)=(x+p\cdot s, y +q\cdot s, s).$$
We will call such a map a {\em Dehn twist}. It is easy to verify that
the extension of $\dd_{(p,q)}$ to $M$ is isotopic to the identity of $M$.
Note that $\dd_{(p,q)}$ is actually not smooth on $T\times\{1\}$, but 
we can consider some smoothing and identify $\dd_{(p,q)}$
to an element of $G(K)$, because the equivalence class is
independent of the smoothing.

\dim{seldom:triv}
We show that $\dd_{(p,q)}$ is non-trivial in $G(K)$ for all $(p,q)\neq(0,0)$. 
Fix the basepoint $a_0=(0,0)\in T$ for the fundamental groups of $T$ and $E(K)$. Then 
$\dd_{(p,q)}$ acts on $\pi_1(E(K),a_0)$ as the conjugation by $i_*(p,q)$, where 
$i:T\to E(K)$ is the inclusion and $(p,q)\in\mz\times\mz=\pi_1(T,a_0)$.
If $\dd_{(p,q)}$ is trivial in $G(K)$, {\em i.e.} it is isotopic to the identity
relatively to $\partial E(K)$, in particular it acts trivially on 
$\pi_1(E(K),a_0)$. This implies that $i_*(p,q)$ is in the centre of $\pi_1(E(K),a_0)$.
Now it follows from hyperbolicity that this centre is trivial and $i_*$ is injective,
whence the conclusion.
\finedim{seldom:triv}

\noindent The proof of Theorem~\ref{many:good} will rely on properties of hyperbolic
manifolds and on the following fact, which we consider to be quite remarkable
(note that the 2-dimensional analogue,
which may be stated quite easily, is false).
Remark that the result applies in particular to Dehn twists.

\begin{prop}\label{collar:trivial}
If $[f]\in G(K)$ and $f$ is supported in the collar $V$ of $\partial U$ then $[f]$
acts trivially on $\eul(E(K),\pp(K^F))$.
\end{prop}

\dim{collar:trivial}
Consider a vector field $v$ on $E(K)$ compatible with $\pp(K^F)$. Since $v$ and $f_*(v)$
differ only on $V$, their difference belongs to the image of $H_1(V;\mz)$ in
$H_1(E(K);\mz)$. So we may as well assume that $E(K)=V$, {\em i.e.} $M$ is
the solid torus $U\cup V$.

By contradiction, let $\xi\in\eul(V,\pp(K^F))$ be such that 
$\alpha(\xi,(\dd_{(p,q)})_*(\xi))$ is non-zero in $H_1(V;\mz)$, so it is given
by $k\cdot[\gamma]$ for some $k\in\mz\setminus\{0\}$ and some simple closed
curve $\gamma$ on $T\times\{1\}\subset\partial V$. Let us now take another
simple closed curve $\delta$ on $T\times\{1\}$ which intersects $\gamma$
transversely at one point. Let us define $N$ as the manifold obtained by
attaching the solid torus to $V$ along $T\times\{1\}$, in such a way that the
meridian of the solid torus gets identified  with $\delta$. Note that $N$ is
again a solid torus and that the homology class of $\gamma$ in
$H_1(N;\mz)\cong\mz$ is a generator. Now we can apply
Proposition~\ref{p:h:formula} to extend $\xi$ to an Euler structure $\xi_N$ on
$N$. Moreover we can extend $f$ to an automorphism $g$ of $N$ which is the
identity on $\partial N=T\times\{0\}$. Now by construction
$\alpha(\xi_N,g_*(\xi_N))$ equals $k\cdot[\gamma]$ in $H_1(N;\mz)\cong\mz$, so
it is non-zero. But $g$ is isotopic to the identity of $N$ relatively to the
boundary of $N$, so we have a contradiction. \finedim{collar:trivial}

For the proof of Theorem~\ref{many:good} we will also need the
following easy fact. 

\begin{lem}\label{power:action}
Let $f$ be an automorphism of $M$ relative to $\partial M$, and consider the
induced automorphisms of $H_1(M;\mz)$ and $\eul(M,\pp)$, both denoted by $f_*$. Then:
$$\alpha(f_*(\xi_0),f_*(\xi_1))=f_*(\alpha(\xi_0,\xi_1)),\qquad
\forall\xi_0,\xi_1\in\eul(M,\pp).$$
\end{lem}

\dim{power:action}
Take representatives of $\xi_0$ and $\xi_1$ such that $\alpha(\xi_0,\xi_1)$ can be
viewed as the anti-parallelism locus. The formula is then obvious.
\finedim{power:action}

\dim{many:good} Consider $[f]\in G(K)$. It follows from the work of Johansson
(see~\cite{hatcher}) that, under the assumption that $E(K)$ is hyperbolic, the
group generated by Dehn twists has finite index in the mapping class group of
$E(K)$ relative to the boundary. More precisely, the quotient group can be
identified to a subgroup of ${\rm Out}(\pi_1(E(K))$, which is finite as a
consequence of Mostow's rigidity. If ${\rm Out}(\pi_1(E(K))$ is trivial then
$[f]$ is equivalent to a Dehn twist, so $f$ acts trivially on
$\eul(E(K),\pp(K^F))$ by Proposition~\ref{collar:trivial}. 

We are left to deal with the case where $H_1(E(K);\mz)$ is torsion-free. By
Johansson's result, there exists an integer $n$ such that $f^n$ acts trivially
on $\eul(E(K),\pp(K^F))$. Consider now $\xi\in\eul(E(K),\pp(K^F))$, and set
$\alpha=\alpha(\xi,f_*(\xi))$.  We must show that  $\alpha=0$. We denote by
$\widehat\alpha$ the image of $\alpha$ in $H_1(M;\mz)$, and by  $\widehat f$
the extension of $f$ to $M$. Since $\widehat f$ is isotopic to the identity,
we  have $\widehat f_*(\widehat\alpha)=\widehat\alpha$. If we take an oriented
1-manifold $a$ representing $\alpha$ and disjoint from $\partial U(K)$, this
means that there exists an  oriented surface $\Sigma$ in $M$ such that
$\partial\Sigma=a\cup(-f(a))$. Up to isotopy we can  assume that $\Sigma$
intersects $\partial U(K)$ transversely in a union of circles. This shows that
$f_*(\alpha)=\alpha+k\cdot\mu$, where $\mu$ is the meridian of $T$. Note that
$f_*(\mu)=\mu$, so for all integers $m$ we have $f_*^m(\alpha)=\alpha+m\cdot
k\cdot\mu$. Now, using Lemma~\ref{power:action}, we have:
	\begin{eqnarray*}
0 & = & \alpha(\xi,f_*^n(\xi)=\sum_{m=0}^{n-1}\alpha(f_*^m(\xi),f_*^{m+1}(\xi))\\
& = & \sum_{m=0}^{n-1} f_*^m(\alpha(\xi,f_*(\xi)))=\sum_{m=0}^{n-1}f_*^m(\alpha)=
\sum_{m=0}^{n-1}(\alpha+m\cdot k\cdot\mu)\\
& = & n\cdot \alpha +{n(n-1)\over 2}\cdot k\cdot \mu.
\end{eqnarray*}
This shows that $2\cdot \alpha+(n-1)\cdot k\cdot \mu$ is a torsion element of $H_1(E(K);\mz)$,
so it is null by assumption. So $(1-n)\cdot k\cdot \mu=2\cdot \alpha$. If we apply $f_*$ to
both sides of this equality we get
$(1-n)\cdot k\cdot f_*(\mu)=2\cdot f_*(\alpha)$. Using the equality again and the relations
$f_*(\mu)=\mu$ and $f_*(\alpha)=\alpha+k\cdot\mu$ we get
$$(1-n)\cdot k\cdot \mu=2\cdot\alpha+2\cdot k\cdot\mu=(1-n)\cdot k\cdot
\mu+2\cdot k\cdot\mu.$$
Therefore $k\cdot\mu$ is a torsion element, and hence null. But
$2\cdot\alpha=(1-n)\cdot k\cdot\mu$, so also $\alpha$ is null.\finedim{many:good}

\paragraph{Rotation number, and more good knots} 
We will show in this section that 
in a homology sphere the rotation number of a pseudo-Legendrian knot can be 
(defined and) expressed 
in terms of Euler structures on the exterior. This will lead us to a 
simple interpretation of torsion as a relative invariant of knots,
and it will allow us to show
that in a homology sphere all knots are good (Theorem~\ref{homology:good}).

To begin, we note that the definition of the rotation number, classically defined in the
contact case, actually extends to the situation we are considering. Since we will need
this definition, we recall it. Let $M$ be a homology sphere, let
$v$ be a field on $M$ and let $K$ be an oriented pseudo-Legendrian knot 
in $(M,v)$. Take a plane field $\eta$ transversal to $v$ and tangent to $K$, 
and a Seifert surface
$S$ for $K$. Up to isotopy of $S$ we can assume that $\eta$ is tangent to
$S$ only at isolated points. Then ${\rm rot}_v(K)$ is the sum of a contribution
for each of these tangency points $p$. Define ${\rm o}(p)$ to be $+1$ if
$\eta_p=T_pS$ and $-1$ if $\eta_p=-T_pS$. If $p\in\partial S=K$ then $p$ contributes
just with ${\rm o}(p)$. If $p\in\interior(S)$ we can consider near $p$ a section of 
$\eta\cap TS$ which vanishes at $p$ only, and denote by ${\rm i}(p)$ its index. Then
$p$ contributes to ${\rm rot}_v(K)$ with ${\rm o}(p)\cdot{\rm i}(p)$.

It is quite easy to see that the resulting number is indeed independent from $\eta$
and $S$. Moreover ${\rm rot}_v(K)$ is unchanged under homotopies of $v$ relative to $K$,
and local modifications away from $K$, so we can actually define
${\rm rot}_\xi(K)$ where $\xi=\xi(v,K)\in\eul(E(K),\pp(K^{(v)})$.

\begin{prop}\label{alpha:rot}
Let $M$ be a homology sphere, let $v$ be a field on $M$ 
and let $K_0$ and $K_1$ be oriented pseudo-Legendrian knots in $(M,v)$. 
Assume that there exists a framed-isotopy $f$ which maps $K_1^{(v)}$ to 
$K_0^{(v)}$. Identify $H_1(E(K_0);\mz)$ to $\mz$ using a meridian. Then:
$${\rm rot}_v(K_1)={\rm rot}_v(K_0)+2\alpha(f_*(\xi(v,K_1)),\xi(v,K_0)).$$
\end{prop}

\dim{alpha:rot}  Let $K:=K_0$, $v_0:=v$ and $v_1:=f_*(v)$. Note that $v_0$ and
$v_1$ coincide along $K$.  Of course ${\rm rot}_v(K_1)={\rm rot}_{v_1}(K)$. We
are left to show that 
$${\rm rot}_{v_1}(K)={\rm rot}_{v_0}(K)+2\alpha(\xi(v_1,K)),\xi(v_0,K)).$$ 
We can now homotope $v_0$ and $v_1$ away from $K$ until they differ only in
the  neighbourhood $W(L)$ of an oriented link $L$, and within this
neighbourhood they differ exactly by a ``Pontrjagin move'', as defined for
instance in~\cite{lnm}. Namely, $v_0$  runs parallel to $L$ in $W(L)$, while
$v_1$ runs opposite to $L$ on $L$ and has non-positive radial component on
$W(L)$ (see below for a picture). Note that $L$ represents
$\alpha(\xi(v_1,K)),\xi(v_0,K))$.

Let us choose now a Seifert surface $S$ for $K$ and a Riemannian metric on $M$,
and define $\eta_i=v_i^\perp$, for $i=0,1$. Since
$\eta_0\ristr{K}=\eta_1\ristr{K}$, the contributions along $K$ to ${\rm
rot}_{v_0}(K)$ and ${\rm rot}_{v_1}(K)$ are the same. Up to isotoping $S$ we
may assume that $L$ is transversal but never orthogonal to $S$. At the points
where $\eta_0$ is tangent to $S$ also $\eta_1$ is tangent to $S$, and the
contributions are the same. So ${\rm rot}_{v_1}(K)-{\rm rot}_{v_0}(K)$ is given
by the sum of the contributions of the tangency points of $\eta_1$ to $S$
within $W(L)$. We will show that each point of $L\cap S$ gives rise to exactly
two tangency points, which both contribute with $+1$ or $-1$ according to the
sign of the intersection of $L$ and $S$ at the point. This will show that ${\rm
rot}_{v_1}(K)-{\rm rot}_{v_0}(K)$ is twice the algebraic intersection of $L$
and $S$. This algebraic intersection is exactly the value of
$[L]=\alpha(\xi(v_1,K)),\xi(v_0,K))$ as a multiple of $[m]$, so the local
analysis at $L\cap S$ will imply the desired conclusion.

\begin{figure}
\centerline{\psfig{file=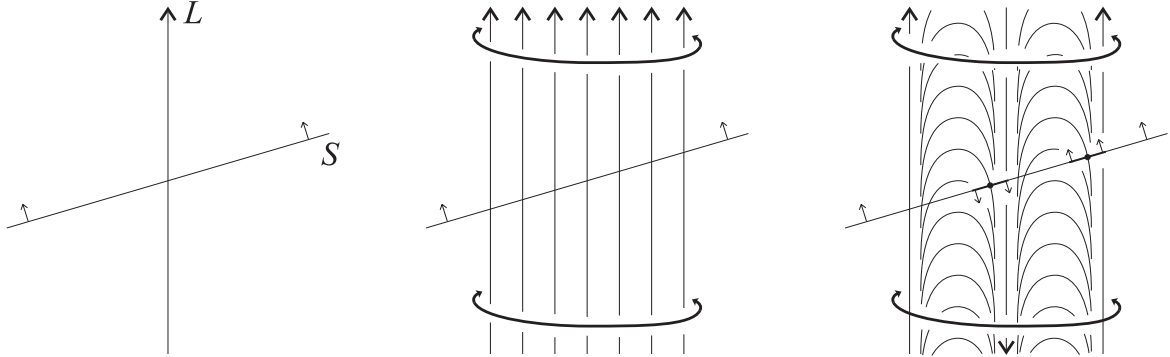,width=15.5cm}}
\caption{\label{rot:alpha:fig}Effect of the Pontrjagin move.}
\end{figure}

For the sake of simplicity we only examine a positive intersection point of $L$
and $S$. This is done in a cross-section in Fig.~\ref{rot:alpha:fig}, which
shows the local effect
of the move. The fields pictured both have a rotational symmetry, suggested in the
figure. The two tangency points which arise are a positive focus (on the right) and
a negative saddle (on the left), 
so the local contribution is indeed $+2$, and the proof is complete.
\finedim{alpha:rot}

\begin{rem}{\em The definition of rotation number and
Proposition~\ref{alpha:rot}
easily extend to the case of manifolds which are not homology spheres, 
by restricting to homologically trivial knots and choosing a relative homology class in the
exterior.}
\end{rem}

\noindent We can now prove that in a homology sphere all knots are good.

\dim{homology:good}
Consider $[f]\in G(K)$, a framing $F$ on $K$ and $\xi\in\eul(E(K),\pp(K^F))$. We must
show that $f_*(\xi)=\xi$. Let $\xi=[v]$ and denote by $\widehat{v}$ the
obvious extension of $v$ to $M$. As above, let $\widehat{f}$ be the extension of $f$ to $M$.
During the proof of Proposition~\ref{alpha:rot} we have shown that
$${\rm rot}_{\widehat{f}_*(\widehat{v})}(K)-{\rm rot}_{\widehat{v}}(K)=
2\alpha(f_*(v),v).$$

But ${\rm rot}_{\widehat{f}_*(\widehat{v})}(K)$ is actually equal to 
${\rm rot}_{\widehat{v}}(K)$, because 
$\widehat{f}$ is the identity near $K$. Therefore
$f_*(v)$ and $v$ differ by a torsion element of
$H_1(E(K);\mz)\cong\mz$, so they are equal. By definition
$f_*(\xi)=[f_*(v)]$ and $\xi=[v]$, and the proof
is complete.
\finedim{homology:good}

Theorems~\ref{many:good} and~\ref{homology:good} provide a partial answer
to the problem of determining which knots are good. The general problem does not appear 
to be straight-forward, and we leave it for further investigation.
We will only show below an example of knot which is not good.

\paragraph{Curls and the winding number}  
We show in this paragraph the relation between the relative invariant 
$\alpha((v_0,K_0),(v_1,K_1))$ of two pseudo-Legendrian knots (when
this invariant is well-defined) and
an analogue of the winding number (the invariant which allows
to distinguish framed-isotopic planar link diagrams which are not equivalent under
the second and third of Reidemeister's moves, see~\cite{trace}). 
Moreover we will give an
example of knot which is not good. The proof of the next result uses the
example of Section~\ref{exa:section}, so it is deferred to
Section~\ref{proofs}.

\begin{prop}\label{wind:sensitive}
Consider a field $v$ on $M$ and a portion of $M$ on which $v$
can be identified to the vertical field in $\mr^3$. Consider oriented knots $K_0$ and
$K_{\pm1}$ which are transversal to $v$ and differ only within the chosen portion of $M$,
as shown in Fig.~\ref{wind:fig}.
\begin{figure}
\centerline{\psfig{file=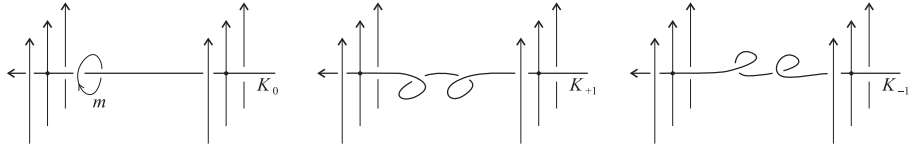,width=12cm}}
\caption{\label{wind:fig}Knots which differ for a positive or a negative double curl.}
\end{figure}
Choose the positive meridian $m$ of $K_0$, as also shown in the figure.
Let $f$ be an isotopy which maps $K_{\pm1}^{(v)}$ to $K_0^{(v)}$ and is supported in
a tubular neighbourhood of $K_0$. Then:
$$\alpha(\xi(v,K_0),f_*(\xi(v,K_{\pm1})))=\pm[m]\in H_1(E(K_0);\mz).$$
\end{prop}

\begin{prop}\label{hyper:curl} Let $(v,K_0)$ be a pseudo-Legendrian pair in
$M$, and denote by $[m]\in H_1(E(K_0);\mz)$ the homology class of the meridian of
$U(K_0)$. Assume either that $K_0^{(v)}$ is good and $[m]\neq 0$ or that
$E(K_0)$ is hyperbolic and $[m]$ has infinite order.
Let $K_1$ be a knot obtained
from $K_0$ as in Fig.~\ref{wind:fig}. Then $(v,K_0)$ and $(v,K_1)$  are not
weakly equivalent. \end{prop}

\dim{hyper:curl} By contradiction, using 
Propositions~\ref{group:obstr} and~\ref{wind:sensitive}, we would get elements
$\xi_0,\xi_1$ of $\eul(E(K_0),\pp(K_0^{(v)})$ such that $\alpha(\xi_0,\xi_1)=[m]$
and $\xi_1=f_*(\xi_0)$ for some $[f]\in G(K_0)$. If $K_0^{(v)}$ is good and $[m]\neq 0$ 
this is a contradiction. Assume now that $E(K_0)$ is hyperbolic and $[m]$ has infinite order.
Since $f_*([m])=[m]$, using 
Lemma~\ref{power:action} we easily see that $\alpha(\xi_0,f_*^k(\xi_0))=k\cdot[m]$ for
all $k$. Proposition~\ref{collar:trivial} and the result of Johansson already used in the proof of
Theorem~\ref{many:good}
now imply that $f^k$ acts trivially on $\eul(E(K_0),\pp(K_0^{(v)})$ for some
$k$, whence the contradiction.\finedim{hyper:curl}

\noindent As an application of Proposition~\ref{wind:sensitive}, 
we can show that there exist knots which are not good. Consider $S^2\times[0,1]$
with vector field parallel to the $[0,1]$ factor. Let $K_0$ be the equator
of $S^2\times\{1/2\}$, and let $K_1$ be obtained from $K_0$ by the
modification described in Fig.~\ref{wind:fig}.
Using Proposition~\ref{wind:sensitive}, if we choose a framed-isotopy
$g$ of $K_1^{(v)}$ onto $K_0^{(v)}$ supported in $U(K_0)$, we have
$$\alpha(\xi(v,K_0),(g\ristr{E(K_1)})_*(\xi(v,K_1)))=[m],$$
where $[m]$ is a generator of $H_1(E(K_0);\mz)\cong\mz$.
On the other hand, $K_1$ is strongly equivalent to $K_0$ in $(M,v)$
(the winding number only exists on $\mr^2$, not on $S^2$). 
So there exists an isotopy $h$ of $K_1^{(v)}$ onto $K_0^{(v)}$ through links
transversal to $v$, and we have
$$\alpha(\xi(v,K_0),(h\ristr{E(K_1)})_*(\xi(v,K_1)))=0.$$
This implies that $(h\compo g^{-1})\ristr{E(K_0)}$ acts non-trivially
on $\xi(v,K_0)\in\eul(E(K_0),\pp(K_0^{(v)}))$.

To conclude our discussion on the relative invariant $\alpha$ between two pseudo-Legendrian
knots, we state now a result proved in~\cite{second:paper}. The consequences we deduce
easily follow from Proposition~\ref{wind:sensitive}.

\begin{prop}\label{curl:diff:prop}
Let $(v_0,K_0)$ and $(v_1,K_1)$ be pseudo-Legendrian in $M$, assume that 
$v_0$ and $v_1$ are homotopic fields, and that $K_0^{(v_0)}$ and $K_1^{(v_1)}$
are isotopic as framed knots. Then $(v_0,K_0)$ and $(v_1,K_1)$ become
weakly equivalent up to a finite number of local moves
$K_0\longrightarrow K_{\pm1}$ as in Fig.~\ref{wind:fig}.
\end{prop}

\begin{cor}\label{rel:wind:cor}
Under the same assumptions, assume also that $K_0^{(v_0)}$ is good, so
$\alpha((K_0,v_0),(K_1,v_1))$ is defined. Then 
$$\alpha((K_0,v_0),(K_1,v_1))=
{\rm w}((K_0,v_0),(K_1,v_1))\cdot[m_0]\in H_1(E(K_0);\mz)$$
where $m_0$ is the meridian of $K_0$ and 
${\rm w}((K_0,v_0),(K_1,v_1))$ is the 
(non-well-defined) algebraic number of moves 
$K_0\longrightarrow K_{\pm1}$ needed to make $K_0$ and $K_1$ weakly equivalent.
\end{cor}

Concerning the statement of the previous corollary, note that 
both ${\rm w}((K_0,v_0),(K_1,v_1))$ and $[m_0]$ depend on the
choice of an orientation on $K_0$, but their product does not.

\begin{cor}
Assume furthermore that $[m_0]$ has infinite order in $H_1(E(K_0);\mz)$. Then 
${\rm w}((K_0,v_0),(K_1,v_1))\in\mz$ is a well-defined integer relative
invariant, which we call the {\em relative winding number}.
\end{cor}

\begin{rem}{\em 
If $M$ is a homology sphere then the local moves of Fig.~\ref{wind:fig}
which modify the winding number also change the rotation number, and 
Corollary~\ref{rel:wind:cor} is consistent with Proposition~\ref{alpha:rot}.}\end{rem}

The next proposition contains in particular Theorem~\ref{informal:good:for:knots} stated
in the introduction.

\begin{prop}\label{no:tors:same}
Under the assumptions of Proposition~\ref{curl:diff:prop}, assume that 
$K_0^{(v_0)}$ is good and that $[m_0]$ has infinite order in 
$H_1(E(K_0);\mz)$. The following facts are pairwise equivalent:
\begin{enumerate}
\item\label{no:wind:point} the relative winding 
number of $(K_0,v_0)$ and $(K_1,v_1)$ vanishes;
\item\label{no:tors:point} all relative torsion 
invariants of $(K_0,v_0)$ and $(K_1,v_1)$ are trivial;
\item\label{equiv:point} $(K_0,v_0)$ and $(K_1,v_1)$ are weakly equivalent.
\end{enumerate}
\end{prop}

\dim{no:tors:same}
Equivalence of (\ref{no:wind:point}) and (\ref{equiv:point})
comes from the previous discussion and from the fact that 
a positive double curl and a negative double curl cancel via
weak equivalence. To show that (\ref{no:wind:point}) and
(\ref{no:tors:point}) are equivalent we only need to consider torsion
with respect to a representation $\varphi:H_1(E(K_0);\mz)\to\Lambda$ such that
$\varphi([m_0])$ has infinite order.
\finedim{no:tors:same}

\begin{cor}\label{immers:cor}
Under the assumptions of Proposition~\ref{curl:diff:prop},
assume that $M$ is a homology sphere. Then the facts (1), (2), and (3) 
of Proposition~\ref{no:tors:same} are also equivalent to the following:
\begin{enumerate}
\item[(4)]\label{no:rot:point} $(K_0,v_0)$ and $(K_1,v_1)$ 
have the same rotation number.
\end{enumerate}
\end{cor}

\dim{immers:cor}
Equivalence of (\ref{no:wind:point}) 
and (\ref{no:rot:point}) comes from the previous discussion and 
Proposition~\ref{wind:sensitive}.\finedim{immers:cor}

Since
in a homology sphere two pseudo-Legendrian knots which are homotopic
through pseudo-Legendrian immersions certainly have the same Maslov index,
the previous corollary seems to suggest that all torsion can capture in a homology
sphere is the homotopy class through immersions.
We believe that it would be interesting to check if also 
for a general manifold $M$, under the assumptions of Corollary~\ref{rel:wind:cor},
homotopy through pseudo-Legendrian immersions implies
${\rm w}((K_0,v_0),(K_1,v_1))\cdot[m_0]=0$. We conclude by informing the reader that
in~\cite{second:paper} we have discussed the extent to which the category of
pseudo-Legendrian knots can be represented by the category of genuine Legendrian
knots in overtwisted contact structures.

\section{An example}\label{exa:section}
\begin{figure}
\centerline{\psfig{file=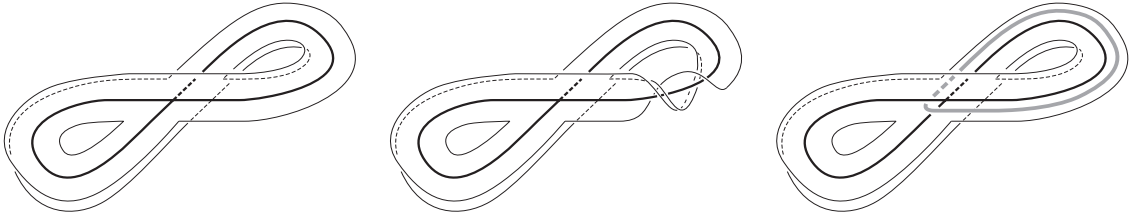,width=15.5cm}}
\caption{\label{abalone}The abalone, and a $\cont^1$ knot on it.}
\end{figure}

Figure~\ref{abalone} shows a neighbourhood of the singular set of the so-called
abalone, a branched standard spine of $S^3$, which we denote by $A$. Note that $A$ has one vertex,
two edges and two regions. The figure on the left is easier to understand, but 
it does not represent the genuine embedding of $A$ in $S^3$, which is instead
shown in the centre (hint: compute linking numbers). On the right we show (using the
easier picture) a $\cont^1$ knot $K$ on $A$. Using the genuine picture one sees that
$K$ is actually trivial in $S^3$, and its framing is $-1$.
So the knot exterior $E(K)$ is actually a solid torus,
with an induced Euler structure $\xi$, and the white annulus $W\subset\partial E(K)$
is a longitudinal one. Let us now take the representation $\varphi:\pi_1(E(K))\to\mz[t^{\pm1}]$
which maps the generator to $t$. It is not hard to see that $H_*^\varphi(E(K),\overline{W})=0$, so we can
compute $\tau^\varphi(E(K),\xi)$. We describe the method to be followed, skipping several details
and all explicit formulae.

\begin{figure}
\centerline{\psfig{file=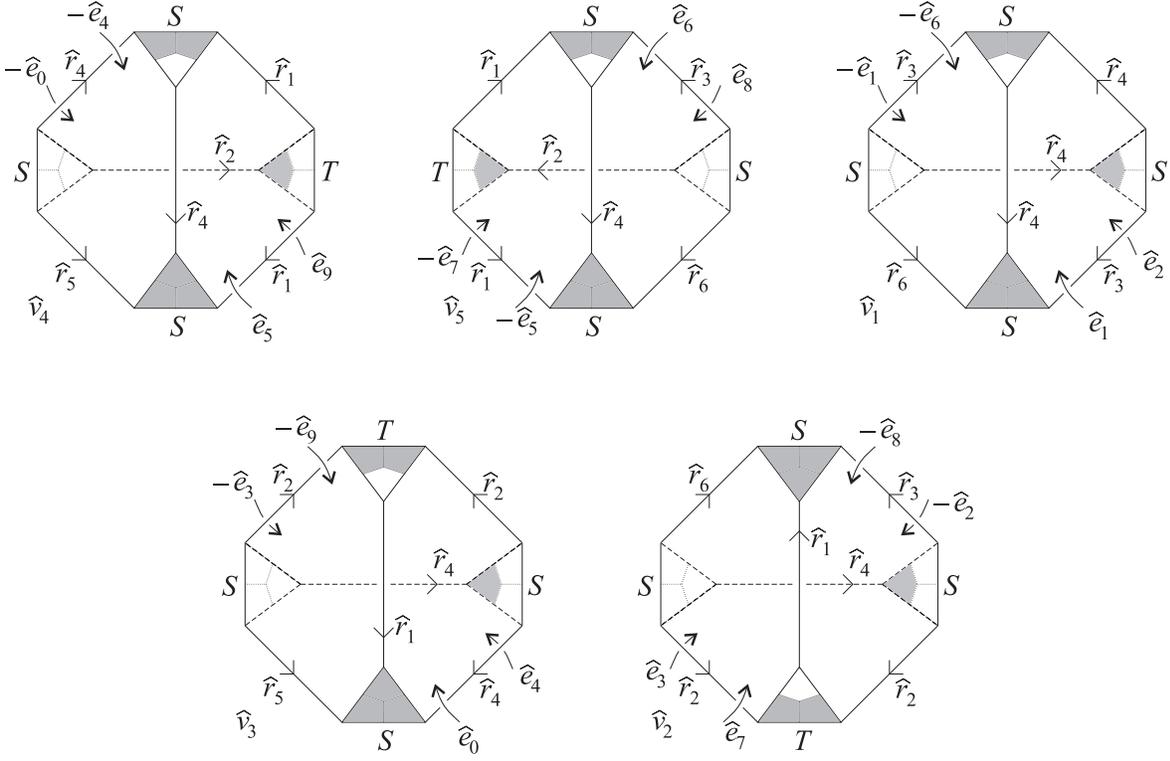,width=15.5cm}}
\caption{\label{dualtria}Truncated ideal triangulation of the knot exterior.}
\end{figure}

We can apply directly the method described in the (partial) proof of 
Theorem~\ref{bounded:surg}, to get a branched standard spine $P$
(in the sense of Theorem~\ref{bounded:surg}) of $E(K)$.
This $P$ is easily recognized to have 5 vertices (denoted $v_1,\dots,v_5$),
10 edges (denoted $e_0,\dots,e_9$) and 6 regions (denoted $r_1,\dots r_6$).
Figure~\ref{dualtria} shows the truncated ideal triangulation dual to $P$.
In the figure the hat denotes duality as usual. We have written
$-\hat e_i$ instead of $\hat e_i$ when $\hat e_i$ lies on 
$\hat v_j$ but the natural
orientation of $\hat e_i$ is not induced by the orientation of $\hat v_j$.
The letters $S$ and $T$ refer to the boundary sphere and torus respectively 
($S$ should actually be collapsed to one point $x_0$, but the picture is
easier to understand before collapse).

Recall that the algebraic complex of which we must compute the
torsion has one generator for each cell in the cellularization 
of $E(K)$ arising from $P$,
excluding the white cells and the tangency circles on the boundary.
From Fig.~\ref{dualtria} we can see how many such cells there will
be in each dimension, namely 3 in dimension 0 ($x_0$ and two vertices
on $T$), 14 in dimension 1 (the $\hat r_i$'s and 8 edges on $T$),
16 in dimension 2 (the $\hat e_i$'s and the 6 black kites on $T$)
and 5 in dimension 3 (the $\hat v_i$'s). We can also easily describe
the combinatorial Euler chain $s'(P)$ which will be used to find
the preferred cell liftings: besides the orbits of the field there are
only one star and one bi-arrow; the support of $s'(P)$ has 3 connected 
components (one spider with 19 legs and head at $x_0$, the star union
the second half of $\hat r_2$, and the bi-arrow union a segment 
contained in $\hat e_3$).

To actually determine the preferred liftings we need an effective
description of the lifting of the cellularization
to the universal cover $\tilde E(K)\to E(K)$.
Since $\pi_1(E(K))=\mz$, each cell $c$ will have liftings $\tilde c^{(n)}$
for $n\in\mz$, where $\tilde c^{(n)}$ is the $n$-th translate of $\tilde c^{(0)}$.
The choice of $\tilde c^{(0)}$ itself is of course arbitrary, but to
understand the cover we must express the $\partial\tilde c^{(0)}$'s in terms of
the other $\tilde d^{(n)}$'s. To do this we start with a lifting $\tilde x_0$
of the basepoint $x_0$ and we lift the other cells one after each other,
taking into account the relations in $\pi_1(E(K)$ and making sure that the union
of cells already lifted is always connected. When a cell $c$ is reached for the first
time, its lifting is chosen arbitrarily and declared to be $\tilde c^{(0)}$, but
its boundary will involve in general $\tilde d^{(n)}$'s with $n\neq 0$.
Once the lifted cellularization is known, it is a simple matter
to determine preferred cell liftings: since the support of $s'(P)$ consists
of 3 spiders, we only need to choose liftings of the 3 heads and then lift
the legs.

Carrying out the computations we have explicitly found the algebraic complex
with coefficients in $\mz[t^{\pm1}]$, and the preferred generators of the
4 moduli appearing. Then, using Maple, we have checked that indeed
the complex is acyclic, and we have computed its torsion as follows:
$$\tau^\varphi(E(K),\xi)=\pm t^{-1}.$$
Note that as an application of Proposition~\ref{wind:sensitive},
 by adding curls, we can 
easily construct a family $\{K_n\}$ of pseudo-Legendrian knots such that
$\tau^\varphi(E(K_n),\xi_n)=\pm t^{n}$.

\section{Main proofs}\label{proofs}

In this section we provide all the proofs which we have
omitted in the rest of the paper. We will always refer to the statements
for the notation.

\dim{p:h:formula}
Let us first recall the classical Poincar\'e-Hopf formula, according to which if $v$ is a
vector field with isolated singularities on a manifold $M$, and $v$ points outwards on
$\partial M$ ({\em i.e.}~$\partial M$ is black), then the sum of the indices of all 
singularities is $\chi(M)$. Assume now that $v$ has isolated singularities and on 
$\partial M$ it is compatible with a pattern $\pp=(W,B,V,C)$. We claim that 
if $\cc$ is a cellularization of $M$ suited to $\pp$ we have:
\begin{equation}\label{ph:formula}
\sum_{x\in{\rm Sing}(v)}\index_x(v)=\chi(M)-\sum_{\sigma\in\cc,\ \sigma\subset W\cup V}
\index(\sigma).
\end{equation}
This formula is enough to prove the statement: if a non-singular field $v$ compatible
with $\pp$ exists then the left-hand side of~\ref{ph:formula} vanishes,
and the right-hand side of~\ref{ph:formula} equals the obstruction of the statement.
On the other hand, if the obstruction vanishes, then one can first consider a singular
field compatible with $\pp$, then group up the singularities in a ball, and remove them.

\begin{figure}
\centerline{\psfig{file=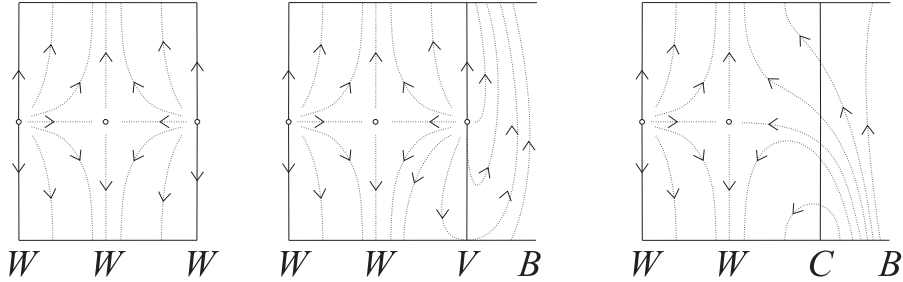,width=12cm}}
\caption{\label{two:obstr}Extension of a field to the collared manifold: dimension 2}
\end{figure}

To prove~\ref{ph:formula} we consider the manifold $M'$ obtained by attaching
a collar $\partial M\times[0,1]$ to $M$ along $\partial M=\partial M\times\{0\}$.
Of course $M'\cong M$. We will now extend $v$ to a field $v'$ on $M'$ with the property that
$v'$ points outwards on $\partial M'$, and in $\partial M\times(0,1)$ the field
$v'$ has exactly one singularity for each cell $\sigma\subset W\cup V$, with index
$\index(\sigma)$. An application of the classical Poincar\'e-Hopf formula then
implies the conclusion. The construction of $v'$ is done cell by cell. We first
show how the construction goes in dimension 2, see Fig.~\ref{two:obstr}.

\begin{figure}
\centerline{\psfig{file=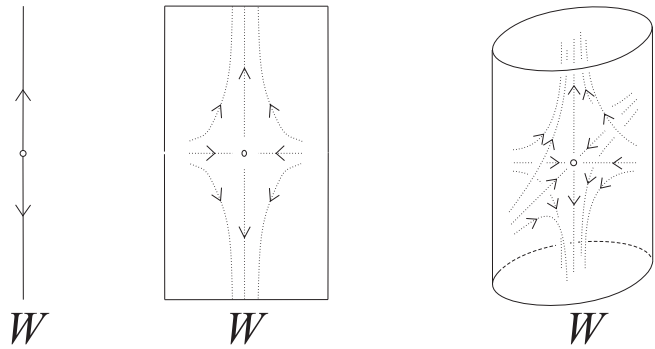,width=9cm}}
\caption{\label{three:obstr:1}Extension of a field to the collared manifold: white
cells in dimension 3}
\end{figure}

\begin{figure}
\centerline{\psfig{file=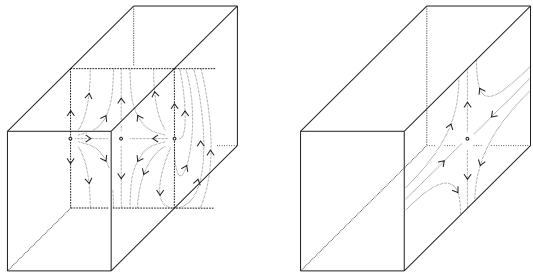,width=7cm}}
\caption{\label{three:obstr:2}Extension of a field to the collared manifold: convex edge in
dimension 3}
\end{figure}

For the 3-dimensional case, it is actually convenient to choose a cellularization $\cc$
of special type. Namely, we assume that $\cc\ristr{\partial M}$ 
consists of rectangles and triangles, each rectangle having exactly one edge on $V\cup C$,
and the union of rectangles covering a neighbourhood of $V\cup C$. We suggest 
in Fig.~\ref{three:obstr:1} how to 
define $v'$ on $\sigma\times[0,1]$ for $\sigma\subset W$ of dimension 0, 1 and 2 respectively.
By the choices we have made the situation near $\partial W$ contains the 2-dimensional
situation as a transversal cross-section, and it is not too difficult to extend $v'$ further
and check that indices of the singularities are as required. As an example, we suggest in
Fig.~\ref{three:obstr:2}
how to do this near a convex edge.\finedim{p:h:formula}

\begin{figure}
\centerline{\psfig{file=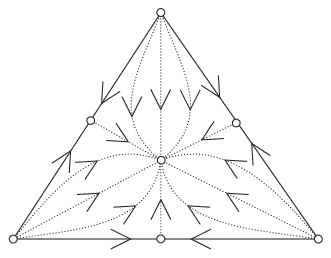,width=4cm}}
\caption{\label{elem:field}The singular field $w_\ss$ on a 2-simplex}
\end{figure}

\dim{reconstruction:statement}
Our proof follows the scheme given by Turaev in~\cite{turaev:Euler}, with some 
technical simplifications and some extra difficulties due to the tangency circles.
We first recall that if $\ss$ is a (smooth) triangulation of a manifold $N$, a (singular)
vector field $w_\ss$ on $N$ can be defined by the requirements 
that: (1) each simplex is a union of orbits; (2)
the singularities are exactly the barycentres of the
simplices; (3) barycentres of higher dimensional simplices are
more attractive that those of lower dimensional simplices.
More precisely, each orbit (asymptotically) goes from a barycentre $p_\sigma$
to a barycentre $p_{\sigma'}$, where $\sigma\subset\sigma'$.
It is automatic that $\index_{p_\sigma}(w_\ss)=\index(\sigma)$.
See Fig.~\ref{elem:field} for a description of $w_\ss$ on a 2-simplex of $\ss$.

Let us consider now a triangulation $\calt$ of $M$, and let us choose a 
representative $z$ of the given $\xi\in\eulc(M,\pp)$ as in
Proposition~\ref{combin:prop:statement}(\ref{bary:sub:point}).
We consider now the manifold $M$ obtained by attaching $\partial M\times[0,\infty)$
to $M$ along $\partial M=\partial M\times\{0\}$. Note that $M'\cong\interior(M)$.
Moreover $\calt$ extends to a ``triangulation'' $\calt'$ of $M'$, where on $M\times(0,\infty)$
we have simplices with exactly one ideal vertex, obtained by taking cones
over the simplices in $\partial M$ and then removing the vertex.
Even if $\calt'$ is not strictly speaking a triangulation, the construction of
$w_{\calt'}$ makes sense, because the missing vertex at infinity would be
a repulsive singularity anyway. We arrange things in such a way that if 
$\sigma\subset\partial M$ then the singularity in $\sigma\times(0,\infty)$
is at height $1$, so it is $p_\sigma\times\{1\}$.

We will define now a smooth function $h:\partial M\to (0,\infty)$ and set $M_h=M\cup\{(x,t)\in\partial M\times[0,\infty):\ t\leq h(x)\}$, in such a way that
$w_{\calt'}$ is non-singular on $\partial M_h$, and, modulo the natural homeomorphism
$M\cong M_h$, it induces on $\partial M_h$ the desired boundary pattern $\pp$.
Later we will show how to use $z$ to remove the singularities of $w_{\calt'}$ on $M_h$.

\begin{figure}
\centerline{\psfig{file=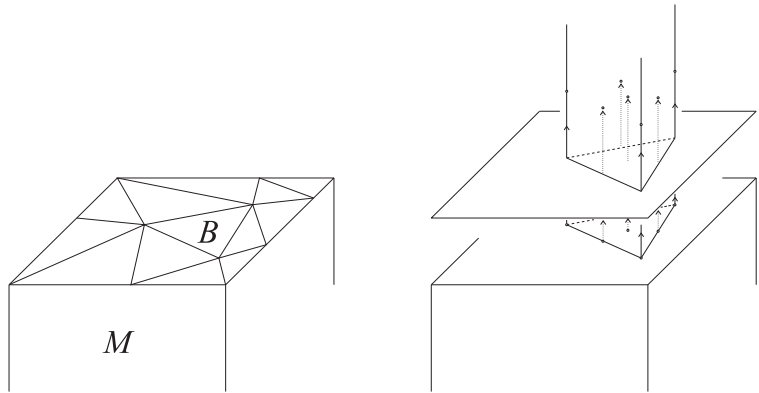,width=10cm}}
\caption{\label{blacktria}Where $h=1/2$ the field points outwards }
\end{figure}

\begin{figure}
\centerline{\psfig{file=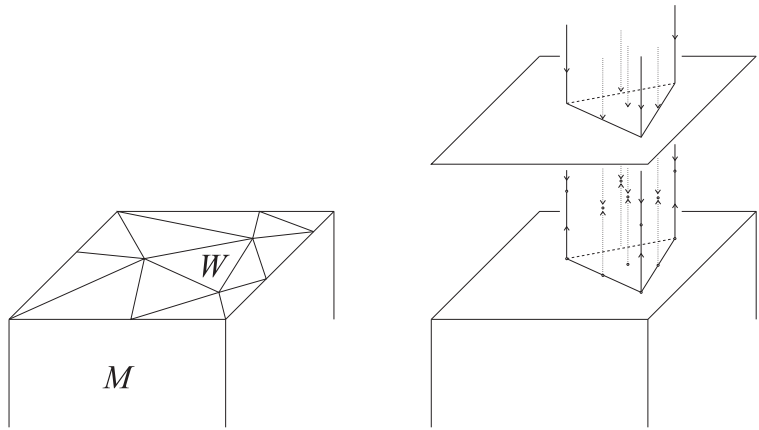,width=10cm}}
\caption{\label{whitetria}Where $h=2$ the field points inwards}
\end{figure}

\begin{figure}
\centerline{\psfig{file=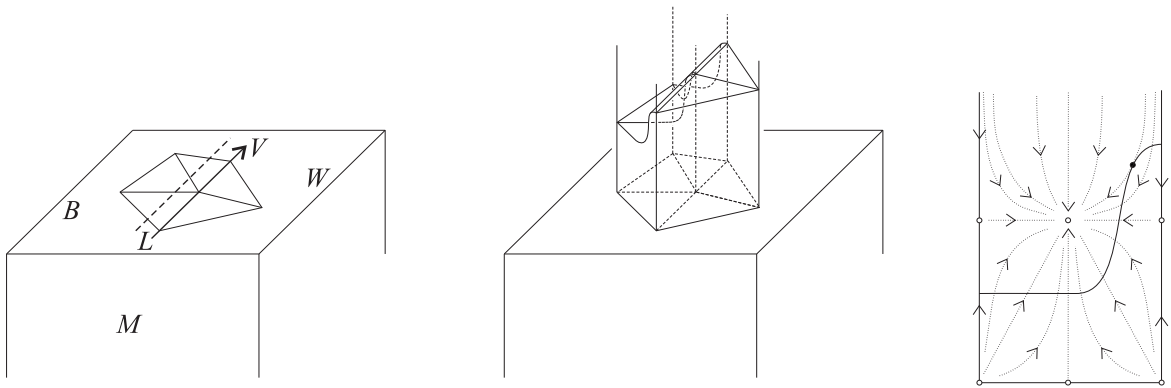,width=15.5cm}}
\caption{\label{convline}On $V$ the field has convex tangency}
\end{figure}

\begin{figure}
\centerline{\psfig{file=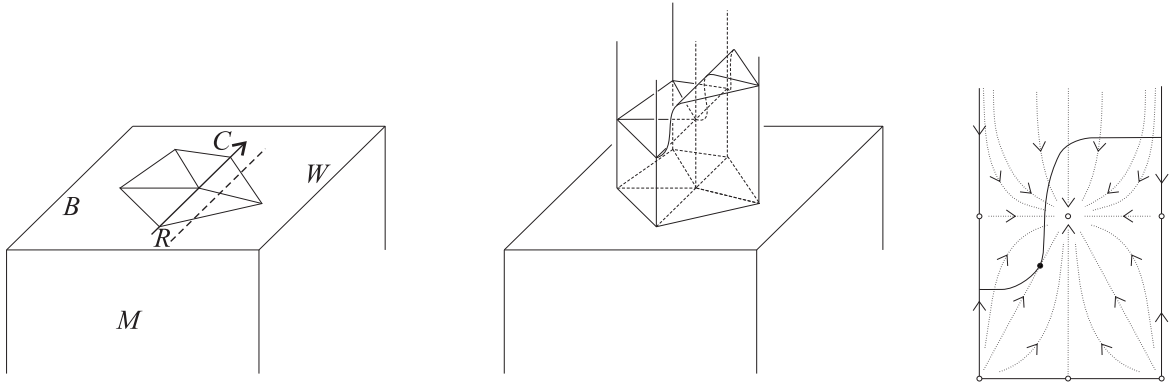,width=15.5cm}}
\caption{\label{concline}On $C$ the field has concave tangency}
\end{figure}

To define the function $h$ we consider a (very thin) left half-collar $L$ of $V$ on 
$\partial M$ and a right half-collar $R$ of $C$. Here ``left'' and ``right'' refer
to the natural orientations of $\partial M$ and of $V\cup C$.
Note that $L\subset B$ and $R\subset W$.
Now we set $h\ristr{B\setminus L}\equiv 1/2$, and 
$h\ristr{W\setminus R}\equiv 2$. Figures~\ref{blacktria} and~\ref{whitetria}
respectively show that away from $V\cup C$ indeed the pattern of $w_{\calt'}$ on
$\partial M_h$ is as required. Now we identify $L$ to $V\times[-1,0]$ and $R$ to
$C\times[0,1]$, and we define $h(x,s)=f(s)$ for $(x,s)\in V\times[-1,0]$
and $h(x,s)=f(s-1)$ for $(x,s)\in C\times[0,1]$, where $f:[-1,0]\to[1/2,2]$ is
a smooth monotonic function with all the derivatives vanishing at $-1$ and $0$.
Instead of describing $f$ explicitly we picture it and show that also near
$V\cup C$ the pattern is as required. This is done near $V$ and $C$ respectively 
in Figg.~\ref{convline}
and~\ref{concline}.
In both pictures we have only considered a special configuration for the 
triangulation on $\partial M$, and we have refrained from picturing the orbits
of the field in the 3-dimensional figure. Instead, we have separately shown
the orbits on the vertical simplices on which the value of $h$ changes.

The conclusion is now exactly as in Turaev's argument, so we only give a sketch.
The chosen representative $z$ of $\xi\in\eulc(M,\pp)$ can be described as an integer linear
combination of orbits of $w_{\calt'}$, which we can describe as segments
$[p_\sigma,p_{\sigma'}]$ with $\sigma\subset\sigma'$. Now we consider the chain
\begin{equation}\label{complete:chain:to:collar}
z'=z-\sum_{\sigma\subset W\cup V}\index(\sigma)\cdot p_\sigma\times[0,1].
\end{equation}
By definition of $h$ we have that $z'$ is a 1-chain in $M_h$, and 
$\partial z'$ consists exactly of the singularities of $w_{\calt'}$ contained in $M_h$, 
each with its index. For each segment $s$ which appear in $z'$
we first modify $w_{\calt'}$ to a field which is ``constant''
on a tube $T$ around $s$, and then we modify the field again
within $T$, in a way which depends on the coefficient of $s$ in $z'$.
The resulting field has the same singularities as $w_{\calt'}$, but 
one checks that these singularities can be removed by 
a further modification supported within small balls centred at the singular points.
We define $\Psi(\xi)$ to be the class in $\euls(M,\pp)$ of this final field.
Turaev's proof that $\Psi$ is indeed well-defined and $H_1(M;\mz)$-equivariant
applies without essential modifications.\finedim{reconstruction:statement}

\begin{rem}\label{good:cell:ok}
{\em In the previous proof we have defined $\Psi$ using triangulations, in order to
apply directly Turaev's technical results (in particular, invariance under subdivision).
However the geometric construction makes sense also for cellularizations $\cc$ more
general than triangulations, the key point being the possibility of defining a field
$w_{\cc}$ satisfying the same properties as the field defined for triangulations.
This is certainly true, for instance, for cellularizations $\cc$ of $M$ induced by
realizations of $M$ by face-pairings on a finite number of polyhedra, assuming
that the projection of each polyhedron to $M$ is smooth.}
\end{rem}

\dim{diagram:commutes}
To help the reader follow the details, we first outline the scheme of the proof:
\begin{enumerate}
\item By identifying $M$ to a collared copy of itself, we choose a representative $z$
of the given $\xi\in\eulc(M,\pp)$ such that the extra terms added to define $\Thetac(\xi)$
cancel with terms already appearing in $z$. (We know {\em a priori} that this happens
at the level of boundaries, but it may well not happen at the level of $1$-chains.)
\item We apply Remark~\ref{good:cell:ok} and choose a cellularization of $M$ in which
it is particularly easy to analyze $\Psi(\xi)$ and $\Psi(\Thetac(\xi))$, both constructed
using the representative $z$ already obtained.
\end{enumerate}

\noindent We consider a cellularization $\cc$ of
$M$ satisfying the same assumptions on $\partial M$ as those considered in the proof of
Proposition~\ref{p:h:formula}, namely $C\cup V$ is surrounded on both sides
by a row of rectangular tiles, and the other tiles are triangular. 
We denote by $\gamma_1,\dots,\gamma_n$ the segments in $C$, oriented as $C$.

Let us consider a representative $z$ relative to $\cc$ of the given $\xi\in\eulc(M,\pp)$.
We construct a new copy $M_1$ of $M$ by attaching $\partial M\times[-1,0]$ to
$M$ along $\partial M=\partial M\times\{-1\}$, and we extend $\cc$ to $\cc_1$ by 
taking the product cellularization on $\partial M\times[-1,0]$. We define a new chain as
\begin{eqnarray*}
z_1 & = & z+\sum_{\sigma\subset B}\index(\sigma)\cdot p_\sigma\times[-1/2,0]
-\sum_{\sigma\subset W\cup V}\index(\sigma)\cdot p_\sigma\times[-1,-1/2]\\
& & + \sum_{j=1}^n\Big(\gamma_j\ristr{[1/2,1]}\times\{-1/2\}-
\gamma_j\ristr{[1/2,1]}\times\{0\}\Big).
\end{eqnarray*}
Note that $z_1$ is an Euler chain in $M_1$ with respect to $\cc_1$. Consider the natural
homeomorphism $f:M\to M_1$ and the class
$$a=\alphac(f_*(\xi),[z_1])\in H_1(M_1;\mz)$$
which may or not be zero. Since the inclusion of $M$ into $M_1$ is an isomorphism
at the $H_1$-level, $a$ can be represented by a $1$-chain in $M$, so $z_1$ can
be replaced by a new Euler chain $z_2$ such that $[z_2]=f_*(\xi)$ and $z_2$
differs from $z_1$ only on $M$.

Renaming $M_1$ by $M$ and $z_2$ by $z$ we have found a representative $z$ of $\xi$
such that $z=z_\theta+\sum_{j=1}^n\gamma_j\ristr{[1/2,1]}$, where $z_\theta$ is a sum
of simplices contained in $B\cup\interior M$. Note that of course $\Thetac(\xi)=[z_\theta]$.
To conclude the proof we need now to analyze $\Psi(\xi)$, constructed using $z$, and
$\Psi(\Thetac(\xi))$, constructed using $[z_\theta]$, and show that
$\Psi(\Thetac(\xi))=\Thetas(\Psi(\xi))$. By construction $\Psi(\xi)$ and 
$\Psi(\Thetac(\xi))$ will only differ near $C$, and we concentrate on one component
of $C$ to show that the difference is exactly (up to homotopy) as in the definition of
$\Thetas$, {\em i.e.}~as in Fig.~\ref{conc:to:conv}.

The difference between $\Psi(\xi)$ and $\Psi(\Thetac(\xi))$ is best visualized on
a cross-section of the form $C\times[0,\infty)$. We leave to the reader to analyze
the complete 3-dimensional pictures. To understand the cross-section, we 
follow the various steps in the proof of Theorem~\ref{reconstruction:statement}.

\begin{figure}
\centerline{\psfig{file=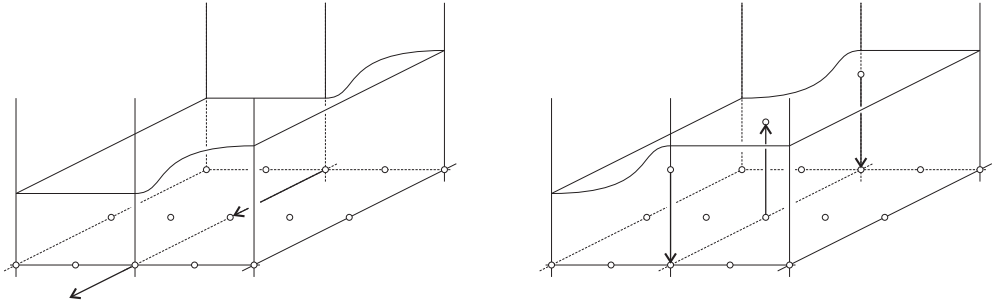,height=4cm}}
\caption{\label{comm:diag:1}Local difference between $z'$ (left) and $z'_\theta$ (right)}
\end{figure}

The first step in the definition of $\Psi(\xi)$ (respectively, $\Psi(\Thetac(\xi))$)
consists in choosing the height function $h$ (respectively, $h_\theta$)
and replacing the chains $z$ (respectively, $z_\theta$)
by a chain $z'$ (respectively, $z'_\theta$)
as in formula~(\ref{complete:chain:to:collar}).
This is done in Fig.~\ref{comm:diag:1}
where only the difference between the chains is shown.

\begin{figure}
\centerline{\psfig{file=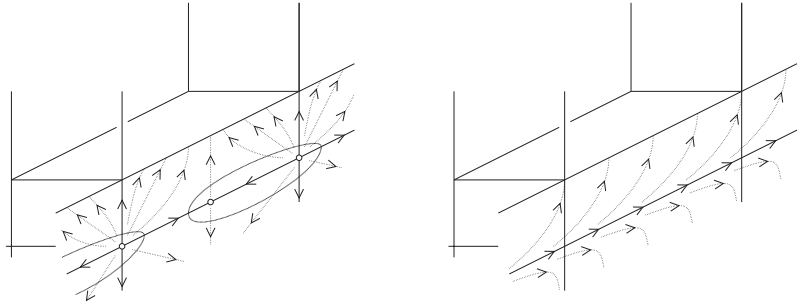,height=4cm}}
\caption{\label{comm:diag:2}Construction of $\Psi(\xi)$ on $C\times[0,\infty)$.
On the left we show $w_\cc$ and the zones where it must be modified.}
\end{figure}

\begin{figure}
\centerline{\psfig{file=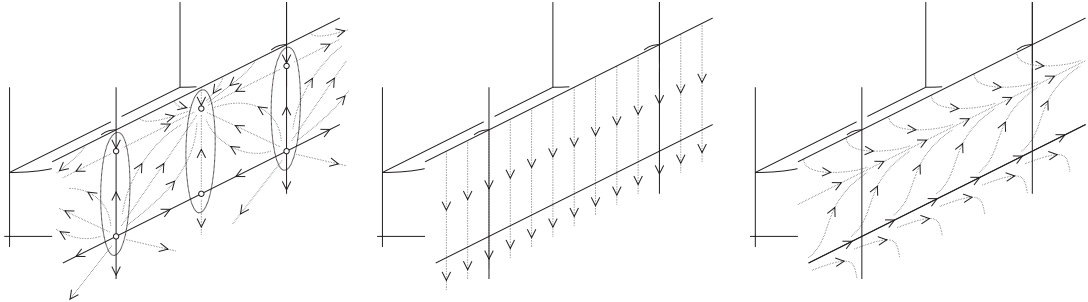,height=4cm}}
\caption{\label{comm:diag:3}Construction of $\Psi(\Thetac(\xi))$ on $C\times[0,\infty)$}
\end{figure}

To conclude we must modify the field $w_\cc$ within a small neighbourhood of 
the support of $z'$ and $z'_\theta$. This is done in Figg.~\ref{comm:diag:2}
and~\ref{comm:diag:3} respectively. The rightmost picture in Fig.~\ref{comm:diag:3} 
is obtained by homotopy on the previous one. The representatives of $\Psi(\xi)$ and
$\Psi(\Thetac(\xi))$ can be compared directly, and indeed they differ by a curve
parallel to $C$ and directed consistently with $C$, so 
$\Psi(\Thetac(\xi))=\Thetas(\Psi(\xi))$.\finedim{diagram:commutes}

We give now the proof omitted in Section~\ref{knots:section}.

\dim{wind:sensitive}
Let us first note that the comparison class which we must show to be $[m]$ is independent
of $f$ by Proposition~\ref{collar:trivial}. We
will give two completely independent (but somewhat sketchy) proofs that this
class is indeed $[m]$.

\begin{figure}
\centerline{\psfig{file=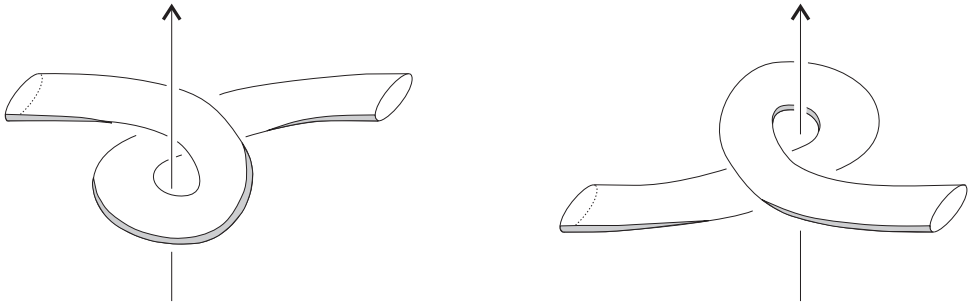,height=4cm}}
\caption{\label{curls}Differently curled tubes in the vertical field.}
\end{figure}

\begin{figure}
\centerline{\psfig{file=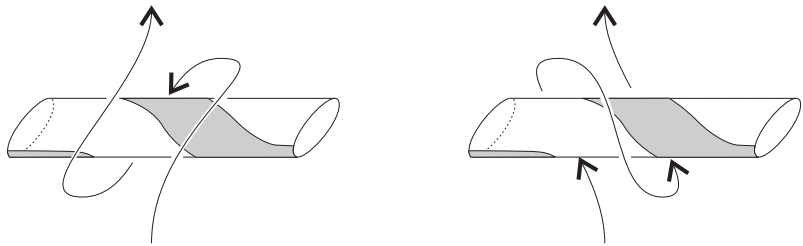,height=3.2cm}}
\caption{\label{nocurls}Straightened curls.}
\end{figure}

For a first proof, instead of comparing a ``straight'' knot with one with two curls, 
we compare two knots with one curl, chosen so that the framing is the same but the 
winding number is different. This is of course equivalent. The two knots are shown in
Fig.~\ref{curls} as thick tubes, together with one specific orbit of the
field they are immersed in. The resulting bicoloration on the boundary
of the tubes is also outlined.
To compare the curls we isotope the bicolorated tubes to the same
straight tube, and we show how the orbit of the field is transformed under
this isotopy. This is done in Fig.~\ref{nocurls}.
Also from this very partial picture it is quite evident that the resulting fields
wind in opposite directions around the tube. A more accurate picture
would show that the difference is actually a meridian of the tube.

Another (indirect) proof goes as follows. Note first that the comparison class which we
must compute certainly is a multiple of $[m]$, say $k\cdot[m]$. Note also that
$k$ is independent of the ambient manifold $(M,v)$. Moreover, by symmetry, we
will have $\alpha(\xi(v,K_0),f_*(\xi(v,K_{-1})))=-k\cdot[m]$ if $K_{-1}$ is
obtained by locally adding a double curl with opposite winding number.

We take now $M$ to be $S^3$, with the field $v$ carried by the abalone $P$
as in Section~\ref{exa:section}, and $K$ to be a trivial knot contained in the 
``smaller'' disc of the $P$. Using either the classical machinery of obstruction
theory or the techniques developed in~\cite{second:paper}, one can see that there
exists another pseudo-Legendrian knot $K'$ in $(S^3,v)$ such that 
$\alpha(\xi(v,K),\xi(v,K'))=[m]$, where by simplicity we are omitting the
framed-isotopies necessary to compare these classes. As already remarked in
Section~\ref{spines:section}, we can assume that $K'$ has a $\cont^1$-projection
on $P$. If one examines $P$ carefully one easily sees that $K'$ can actually be slid
over $P$ to lie again in the small disc of $P$. Now $K'$ is a planar projection of the
trivial knot, so through Reidemeister moves of types II and III,
which correspond to isotopies through knots transversal to $v$, it can be transformed
into a projection which differs from the trivial one only for a finite (even) number
of curls. Summing up, we have a knot $K'$ such that $\alpha(\xi(v,K),\xi(v,K'))=[m]$
and $K'$ differs from $K$ only for a finite number of transformations of the form
$K\mapsto K_1$ or $K\mapsto K_{-1}$. This shows that $[m]$ is a multiple of 
$k\cdot[m]$, so $k=\pm1$.\finedim{wind:sensitive}

We conclude the paper by establishing
the only statement given in Section~\ref{spines:section} and
not proved there. As above, we do not recall all the notation.

\dim{spider:structure:teo}
We fix $P$ and set $s''=s''(P)$, $\hatv=\hatv(P)$.
Using Remark~\ref{good:cell:ok} we see that the construction of
$\Psi([s''])$ explained in the proof of Theorem~\ref{reconstruction:statement}
can be directly applied to the cellularization $\hattt=\hattt(P)$ of $\hatM$.
Recall that this construction requires identifying $\hatM$ to a collared copy
of itself, and extending $s''$ to a chain $s'''$ whose boundary consists precisely
of the singularities of a field $w$. A representative of $\Psi([s''])$ is then
obtained by applying to $w$ a certain desingularization procedure.
This desingularization is supported in a
neighborhood of $s'''$, and one can easily check that
the connected components of the support of $s'''$ (denoted henceforth by $S$)
are actually contractible.
Therefore, {\em any} desingularization of $w$ supported in a
neighbourhood of $s'''$ will give a representative of $\hatv$.
We will prove the desired formula $\Psi([s''])=[\hatv]$
by exhibiting one such desingularization which is nowhere antipodal to $\hatv$.
In our argument we will always neglect the contraction of $\stwotriv$ which maps
$M$ onto $\hatM$. (The desired formula actually holds at the level of $M$, and it
easily implies the formula for $\hatM$.) 

\begin{figure}
\centerline{\psfig{file=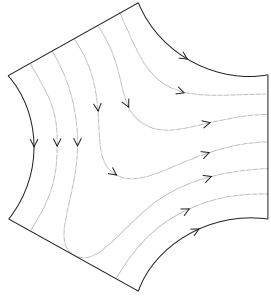,width=3.5cm}}
\caption{\label{hexagon:1}The field $\hatv$ on a hexagon}
\end{figure}

\begin{figure}
\centerline{\psfig{file=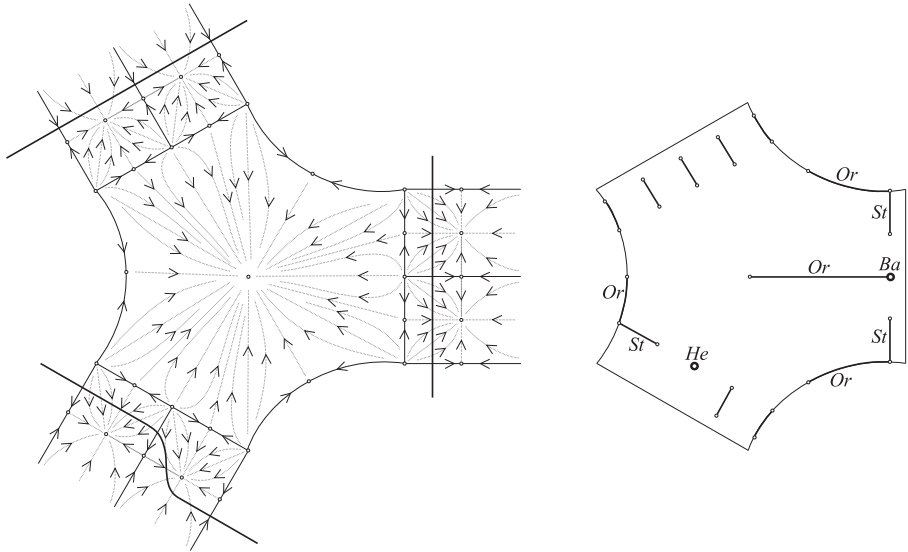,width=12cm}}
\caption{\label{hexagon:2}The field $w$ and the trace of $S$ on a hexagon}
\end{figure}

By the above observations, the following claims easily imply the conclusion of the proof:
\begin{enumerate}
\item\label{claim:antip} The set of points where $w$ is antipodal to 
$\hatv$ is contained in $S$.
\item\label{claim:desing} If $S_0$ is a component of
$S$ then $w$ can be desingularized within a neighbourhood of
$S_0$ to a field which is not antipodal to $\hatv$ in the neighbourhood.
\end{enumerate}
\noindent We prove claim~\ref{claim:antip} by first noting that the cells dual
to those of $P$ are unions of orbits of both $w$ and $\hatv$. So we can analyze cells
separately. We do this explicitly only for 2-dimensional cells, leaving to 
the reader the other cases. In Fig.~\ref{hexagon:1}
we describe $\hatv$. In the left-hand side of Fig.~\ref{hexagon:2}
we describe $w$ on the collared hexagon. In the right-hand side of the same figure
we only show the singularities of $w$ on the renormalized hexagon,
and the intersection of $S$ with the hexagon. In this figure the 7 short segments come
from $s'''-s''$; the other bits of $S$ have been labeled by `Or', `St', `Ba' or `He'
to indicate that they come from
orbits of $\hatv$, stars, bi-arrows or half-edges.

\begin{figure}
\centerline{\psfig{file=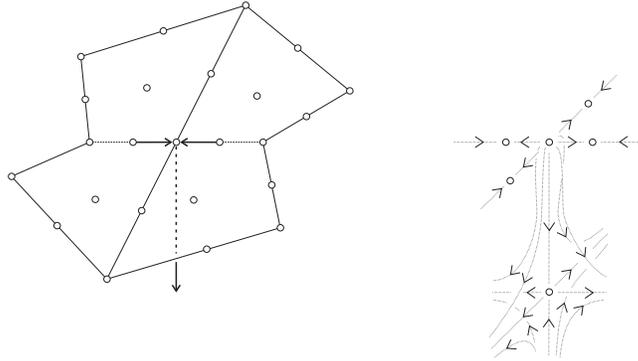,width=8.5cm}}
\caption{\label{fix:biarrow}An enhanced bi-arrow and the field $w$ near it}
\end{figure}

This proves claim~\ref{claim:antip}. Comparing Fig.~\ref{hexagon:2}
with Fig.~\ref{hexagon:1}, and carrying out the same analysis for
3-cells, one actually shows also claim~\ref{claim:desing} for 
components $S_0$ coming from $s'''-s''$. Components of $S$ other than these can
be described in one of the following ways:
\begin{itemize}
\item[(a)] an orbit of $\hatv$ emanating from a vertex of $P$;
\item[(b)] a half-edge of $C$;
\item[(c)] a bi-arrow together with an orbit of $\hatv$ emanating from
the centre of an edge of $P$ and reaching the centre of the bi-arrow;
\item[(d)] a star together with an orbit of $\hatv$ emanating from
the centre of a disc of $P$ and reaching the centre of the star.
\end{itemize}
\noindent All cases can be treated with the same method, we only do case (c). 
Figure~\ref{fix:biarrow}
shows the component placed so that $\hatv$ can be thought of as the 
constant vertical field pointing upwards, and the field $w$ near the component.
The conclusion easily follows.
\finedim{spider:structure:teo}

\vspace{1cm}

\hspace{8cm} benedett@dm.unipi.it

\hspace{8cm} petronio@dm.unipi.it

\hspace{8cm} Dipartimento di Matematica

\hspace{8cm} Via F.~Buonarroti, 2

\hspace{8cm} I-56127, PISA (Italy)


\begin{thebibliography}{99}




\bibitem[1]{manuscripta} {\sc R.~Benedetti, C.~Petronio}, {\it A finite
graphic calculus  for $3$-manifolds}, Manuscripta Math. {\bf 88} (1995),
291-310.

\bibitem[2]{lnm}  {\sc R.~Benedetti, C.~Petronio}, ``Branched Standard
Spines of 3-Manifolds'', Lecture Notes in Math. n. 1653, Springer-Verlag,
Berlin-Heidelberg-New York, 1997.

\bibitem[3]{contspin} {\sc R.~Benedetti, C.~Petronio},  {\it Branched spines
and contact structures on $3$-manifolds}, Ann. Mat. Pura Appl. (to appear).

\bibitem[4]{second:paper} {\sc R.~Benedetti, C.~Petronio},  
{\it Combed 3-manifolds with concave boundary,
framed links, and pseudo-Legendrian links},
preprint January 2000, {\tt math.GT/0001162 }, submitted.

\bibitem[5]{third:paper}  {\sc R.~Benedetti, C.~Petronio},  
{\it Reidemeister torsion via branched spines}, (in preparation).

\bibitem[6]{fourth:paper} {\sc R.~Benedetti, C.~Petronio},  {\it
Euler structures with generic boundary and torsion of
Reidemeister type}, (in preparation).

\bibitem[8]{casler} {\sc B.~G.~Casler}, {\it An imbedding theorem for
connected $3$-manifolds with boundary}, Proc. Amer. Math. Soc. {\bf 16} (1965),
559-566.

\bibitem[9]{christy} {\sc J.~Christy}, {\it Branched surfaces and attractors
I}, Trans. Amer. Math. Soc. {\bf 336} (1993), 759-784.

\bibitem[11]{fist} {\sc R.~Fintushel, R. Stern}, {\it Knots, links,
and $4$-manifolds}, Invent. Math. {\bf 134} (1998), 363-400.

\bibitem[12]{gr:1} {\sc D.~Gillman, D.~Rolfsen}, {\it The Zeeman conjecture
is equivalent to the Poincar\'e conjecture}, Topology {\bf 22} (1983), 315-323.

\bibitem[13]{hatcher} {\sc A.~Hatcher, D.~McCullough}, {\it Finiteness
of classifying spaces of relative diffeomorphism groups of $3$-manifolds}, 
Geometry and Topology {\bf 1} (1997), 91-109.

\bibitem[14]{ishii} {\sc I.~Ishii}, {\it Moves for flow-spines and
topological invariants of $3$-manifolds}, Tokyo J. Math. {\bf 15} (1992),
297-312.

\bibitem[15]{lemuoh}{\sc T.~T.~Q.~Le, J.~Murakami, T.~Ohtsuki},
{\it On a universal perturbative invariant of $3$-manifolds},
Topology {\bf 37} (1998), 359-374.

\bibitem[16]{lescop}{\sc C.~Lescop},``Global Surgery Formula for the Casson-Walker invariant'',
Ann. of Math. Studies n. 140, Princeton University Press, Princeton, NJ, 1996.

\bibitem[18]{mafo} {\sc S.~V.~Matveev, A.~T.~Fomenko},
{\it Constant energy surfaces of Hamiltonian systems, enumeration
of three-dimensional manifolds in increasing order of complexity,
and computation of volumes of closed hyperbolic manifolds},
Russ. Math. Surv. {\bf 43} (1988), 3-25.

\bibitem[19]{matv:mossa} {\sc S.~V.~Matveev}, {\it Transformations of special
spines and the Zeeman conjecture}, Math. USSR-Izv. {\bf 31} (1988), 423-434.

\bibitem[20]{matv:new}
{\sc S.~V.~Matveev}, ``Algorithmic Methods in
3-Manifold Topology'', in preparation.

\bibitem[21]{meng:taub}{\sc G.~Meng, C.~H.~Taubes},
{\it $\underline{SW}\,=\,$Milnor torsion}, Math. Res. Lett. {\bf 3} (1996), 661-674

\bibitem[22]{milnor}{\sc J.~Milnor}, {\it Whitehead torsion}, Bull. Amer. Math. Soc
{\bf 72} (1966), 358-426.

\bibitem[23]{tesi} {\sc C.~Petronio}, ``Standard Spines and 3-Manifolds'', Scuola
Normale Superiore, Pisa, 1995.

\bibitem[24]{piergallini} {\sc R.~Piergallini}, {\it Standard moves for  
standard polyhedra and spines}, Rendiconti Circ. Mat. Palermo {\bf 37}, suppl.
18 (1988), 391-414.

\bibitem[25]{porti} {\sc J.~Porti},  ``Torsion de Reidemeister pour les
Vari\'et\'es Hyperboliques'', Memoirs n. 612, Amer. Math. Soc., Providence, RI,
1997.

\bibitem[26]{trace} {\sc B.~Trace}, {\it On the Reidemeister
moves of a classical knot}, Proc. Amer. Math. Soc. {\bf 89} (1983), 722-724.

\bibitem[27]{turaev:Reidemeister} {\sc V.~G.~Turaev}, {\it Reidemeister
torsion in knot theory}, Russ. Math. Surv. {\bf 41} (1986), 119-182.

\bibitem[28]{turaev:Euler} {\sc V.~G.~Turaev}, {\it Euler structures, nonsingular
vector fields, and torsion of Reidemeister type}, Math. USSR-Izv. {\bf 34}
(1990), 627-662.

\bibitem[29]{turaev:spinc} {\sc V.~G.~Turaev}, {\it Torsion invariants of
Spin$^{\it c}$-structures on $3$-manifolds}, Math. Res. Lett. {\bf 4} (1997), 679-695.

\bibitem[30]{turaev:nuovo} {\sc V.~G.~Turaev}, {\it A combinatorial
formulation for Seiberg-Witten invariants of $3$-manifolds},
Math. Res. Lett. {\bf 5} (1998), 583-598.



\end{thebibliography}
\end{document}